\documentclass[11pt]{article}

\usepackage{times}
\usepackage{float}

\usepackage{amsfonts,amstext,amsmath,amssymb,amsthm}

\usepackage{chngpage}
\usepackage{rotating}
\usepackage[overload]{textcase}
\usepackage{graphicx}
\usepackage{epstopdf} 

\usepackage[hang]{subfigure}
\usepackage[hang]{caption2}

\long\def\symbolfootnote[#1]#2{\begingroup%
\def\thefootnote{\fnsymbol{footnote}}\footnote[#1]{#2}\endgroup}

\usepackage{titlesec} 

\titleformat{\section}{\large\bfseries\uppercase}{\thesection.}{1em}{}
\titlespacing*{\section}{0pt}{*3}{*2}
\titleformat{\subsection}{\normalfont\bfseries}{\thesubsection.}{.5em}{}
\titlespacing*{\subsection}{0pt}{*3}{*2}
\titleformat{\subsubsection}{\normalfont\bfseries}{\thesubsubsection.}{.5em}{}
\titlespacing*{\subsubsection}{0pt}{*3}{*2}

\usepackage{nameref}

\usepackage[%
    colorlinks=true,%
    bookmarksnumbered=true,%
    bookmarksopen=true,%
    citecolor=blue,%
    urlcolor=blue,%
    unicode=true,           
    breaklinks=true,        
    pdftitle={Nearly Optimal Change-Point Detection with an Application to Cybersecurity},%
    pdfauthor={Aleksey S. Polunchenko, Alexander G. Tartakovsky and Nitis Mukhopadhyay},%
]{hyperref}

\numberwithin{equation}{section}

\usepackage{paralist}

\usepackage[sort&compress]{natbib}

\newcommand{\ignore}[1]{}

\renewcommand{\Pr}{\mathbb{P}} 
\DeclareMathOperator{\EV}{\mathbb{E}} 

\DeclareMathOperator{\LR}{\Lambda}

\DeclareMathOperator{\ADD}{ADD}
\DeclareMathOperator{\ARL}{ARL}

\DeclareMathOperator{\IADD}{IADD}

\DeclareMathOperator*{\arginf}{arg\,inf}

\DeclareMathOperator*{\esssup}{ess\,sup}

\newcommand{\tinyinfty}{\scriptscriptstyle\infty}

\newcommand{\T}{T}

\renewcommand{\le}{\leqslant} 
\renewcommand{\ge}{\geqslant}

\theoremstyle{plain} 

\theoremstyle{remark}

\newtheorem*{remark*}{Remark}

\theoremstyle{definition} 

\newtheorem*{definition*}{Definition}

\usepackage{color}
\definecolor{deprecated}{rgb}{0.5,0.5,0.5}

\graphicspath{{gfx/}}

\begin{document}

\title{\textbf{\Large Nearly Optimal Change-Point Detection with an Application to Cybersecurity}}

\date{}
\maketitle

\begin{center}
\null\vskip-2cm\author{
\textbf{\large Aleksey\ S.\ Polunchenko\ and\ Alexander\ G.\ Tartakovsky}\\
Department of Mathematics, University of Southern California, Los Angeles, California, USA
\vskip0.2cm
\textbf{\large Nitis Mukhopadhyay}\\
Department of Statistics, University of Connecticut, Storrs, Connecticut, USA
}
\end{center}
%
%
%
%
\symbolfootnote[0]{\normalsize\hspace{-0.6cm}Address correspondence to A.\ G.\ Tartakovsky, Department of Mathematics, University of Southern California, KAP 108, Los Angeles, CA 90089-2532, USA; Tel: +1 (213) 740-2450, Fax: +1 (213) 740-2424; E-mail:~\href{mailto:tartakov@math.usc.edu}{tartakov@math.usc.edu}.}\\
%
%
{\small\noindent\textbf{Abstract:} We address the sequential change-point detection problem for the Gaussian model where baseline distribution is Gaussian with variance $\sigma^2$ and mean $\mu$ such that $\sigma^2=a\mu$, where $a>0$ is a known constant; the change is in $\mu$ from one known value to another. First, we carry out a comparative performance analysis of four detection procedures: the CUSUM procedure, the Shiryaev--Roberts (SR) procedure, and two its modifications  -- the Shiryaev--Roberts--Pollak and Shiryaev--Roberts--$r$ procedures. The performance is benchmarked via Pollak's maximal average delay to detection and Shiryaev's stationary average delay to detection, each subject to a fixed average run length to false alarm. The analysis shows that in practically interesting cases the accuracy of asymptotic approximations is ``reasonable'' to ``excellent''. We also consider an application of change-point detection to cybersecurity -- for rapid anomaly detection in computer networks. Using real network data we show that statistically traffic's intensity can be well-described by the proposed Gaussian model with $\sigma^2=a\mu$ instead of the traditional Poisson model, which requires $\sigma^2=\mu$. By successively devising the SR and CUSUM procedures to ``catch'' a low-contrast network anomaly (caused by an ICMP reflector attack), we then show that the SR rule is quicker. We conclude that the SR procedure is a better cyber ``watch dog'' than the popular CUSUM procedure.
}
\\ \\
%
%
\\ \\
%
%
%
%

\section{Introduction}
\label{sec:intro}

Sequential change-point detection is concerned with the design and analysis of techniques for fastest (on-line) detection of a change in the state of a process, subject to a tolerable limit on the risk of committing a false detection. The basic iid version of the problem considers a series of independent random observations, $X_1,X_2,\ldots$, which initially follow a common, known pdf $f(x)$, but subsequent to some unknown time index $\nu$ all adhere to a common pdf $g(x)\not\equiv f(x)$, also known; the unknown time index, $\nu$, is referred to as the change-point, or the point of ``disorder''. The objective is to decide after each new observation whether the observations' common pdf is currently $f(x)$, and either continue taking more data, if the decision is positive, or stop and trigger an ``alarm'', if the decision is negative. At each stage, the decision is made based solely on the data observed up to that point. The problem is that the change is to be detected with as few observations as possible past the true change-point, which must be balanced against the risk of sounding a ``false alarm'' -- stopping prematurely as a result of an erroneously made conclusion that the change did occur, while, in fact, it never did. That is, the essence of the problem is to reach a tradeoff between the loss associated with the detection delay and the loss associated with false alarms. A good sequential detection procedure is expected to minimize the average detection delay, subject to a constraint on the false alarm risk.

The particular change-point scenario we study in this work assumes that
\begin{align}\label{eq:Gauss2Gauss-model}
f(x)
&=
\frac{1}{\sqrt{2\pi a\mu}}\exp\left\{-\frac{(x-\mu)^2}{2a\mu}\right\}
\quad\text{and}\quad g(x)=\frac{1}{\sqrt{2\pi a\theta}}\exp\left\{-\frac{(x-\theta)^2}{2a\theta}\right\},
\end{align}
where $0<a<\infty$, $0<\mu\neq\theta<\infty$, and $-\infty<x<\infty$. We will refer to this model as the $\mathcal{N}(\mu,a\mu)$-to-$\mathcal{N}(\theta,a\theta)$ model. This model can be motivated by several interesting application areas.

First, it arises in the area of Statistical Process Control (SPC) in a sequential estimation context.
Suppose a process working in the background on day $i\ge1$ produces iid data $\{X_{i,n}\}_{n\ge1}$ with unknown mean, $\mu$, and unknown variance, $\sigma^{2}$. Per an internal SPC protocol one would like to estimate $\mu$ by constructing a confidence interval for it of prescribed width $2d$, $d>0$, and prescribed confidence level $1-\alpha$, $0<\alpha<1$. If $\sigma^2$ were known, then the required asymptotically (as $d\to0$) optimal (fixed) sample size would be $N_{\alpha,d}^*=\lceil z_{\alpha/2}^{2}\sigma^{2}/d^{2}\rceil$, where $\lceil x\rceil$ is the smallest integer not less than $x$. Specifically, we would have $\lim_{d\to0}\bigl[d^2 N_{\alpha,d}^*/z_{\alpha/2}^2\sigma^2\bigr]=1$, and therefore, $\lim_{d\to0}\Pr(|\bar{X}_{N_{\alpha,d}^*}-\mu|\le d)=1-\alpha$, where $\bar{X}_n=(1/n)\sum_{i=1}^n X_i$ is the sample mean.
However, since $\sigma^{2}$ is unknown,~\cite{Chow+Robbins:AMS1965} proposed the following purely sequential strategy. Let
\begin{align*}
N_{\alpha,d}(i)
&=
\inf\left\{n\ge m_{i}(\ge 2)\colon n\ge z_{\alpha/2}^{2}\frac{S_{i,n}^{2}}{d^{2}}\right\},
\end{align*}
where $\{S_{i,n}^{2}\}_{n\ge2}$ is the sample variance computed from $\{X_{i,n}\}_{n\ge1}$ on day $i\ge1$, and $m_{i}$ is the pilot sample size on day $i\ge1$. At the end of the $i$-th day, we get to observe $N_{\alpha,d}^*(i)$, and thus obtain an estimate of $N_{\alpha,d}^*$. Note that $N_{\alpha,d}(i)$, $i\ge1$, is finite w.p. 1. Furthermore, under the sole assumption that $0<\sigma^2<\infty$,~\cite{Chow+Robbins:AMS1965} show that $N_{\alpha,d}(i)$ is asymptotically (as $d\to0$) efficient and consistent. The same result when the data are Gaussian was also established by~\cite{Anscombe:JRSS53,Anscombe:CPhS62}. Now, by~\citet[Theorem~2.4.1, p.~41]{Mukhopadhyay+Solanky:Book94}, first established by~\cite{Ghosh+Mukhopadhyay:India75}, as $d\to0$, the distribution of $N_{\alpha,d}(i)$ converges to Gaussian with mean $N_{\alpha,d}^*$ and variance $a N_{\alpha,d}^*$, where $a>0$ is foundable explicitly; for instance, $a=2$ when the data are Gaussian. One may additionally refer to \citet[Exercise~2.7.4, p.~66]{Ghosh+etal:Book97} and~\cite{Mukhopadhyay+deSilva:Book09}.

We note that when sampling in multi-steps or under some sort of acceleration parameter $0<\rho<1$, this ``$a$'' is a function of $\rho$, and captures how much sampling is done purely sequentially before it is augmented by batch sampling. See, e.g.,~\citet[Chapter~6]{Ghosh+etal:Book97}.

The main concern is $\sigma^{2}$, that is, by trying to hold $N_{\alpha,d}^*$ within reason, one may want to detect changes in $N_{\alpha,d}^*$ over days. Put otherwise, what if we suddenly see a trend of ``larger or smaller than usual'' values of $N_{\alpha,d}(i)$, $i\ge1$, i.e., estimates of $N_{\alpha,d}^*$? As $d\to0$, this becomes precisely model~\eqref{eq:Gauss2Gauss-model}. This discussion remains unchanged even when the mean $\mu=\mu_{i}$, $i\ge1$, is not the same over days.

Second, model~\eqref{eq:Gauss2Gauss-model} can find application in telecommunications. Let $X_n$ and $Y_n$ be the number of calls made to location 1 and 2, respectively, on day $n\ge1$. Assume that $\{X_n\}_{n\ge1}$ and $\{Y_n\}_{n\ge1}$ are each iid Poisson with rate $\lambda>0$. At location 1, the number of calls, $X_n$ is recorded internally during the whole day's operation between 8:00 AM to 6:00 PM. At location 2, however, the number of calls can be monitored only during half-the-day, between 8:00 AM to 1:00 PM, e.g., for lack of appropriate fund or staff. How should one model the aggregate number of calls, i.e., the number of calls from lines 1 and 2 combined? Since $X_n+Y_n$ cannot be observed due to lack of PM staffing, one may instead consider $U_n=X_n+Y_n/2$. This is reasonable provided the calls are divided approximately equally between the morning and afternoon shifts. For moderately large $\lambda$, however, $U_n$ will be approximately Gaussian with mean $3\lambda/2$ and variance $5\lambda/4$. Thus, at the end of day $i\ge1$ one respectively has independent observations $U_{1},U_{2},\ldots$ recorded from a $\mathcal{N}(\mu,a\mu)$ distribution with $\mu=3\lambda/2$, $a=5/6$.

Third, the possibility of approximating the Poisson distribution with Gaussian makes model~\eqref{eq:Gauss2Gauss-model} of use in the area of cybersecurity, specifically in the field of rapid volume-type anomaly detection in computer networks traffic. A volume-type traffic anomaly can be  caused by many reasons, e.g., by an attack (intrusion), a virus, a flash crowd, a power outage, or misconfigured network equipment. In all cases it is associated with a sudden change (usually, from lower to higher) in the traffic's volume characteristics; the latter may be defined as, e.g., the traffic's intensity measured via the packet rate -- the number of network packets transmitted through the link per time unit. This number randomly varies with time, and according to latest findings, its behavior can be well described by a Poisson process; see, e.g.,~\cite{Cao+etal:NEC02},~\cite{Karagiannis+etal:IEEE-INFOCOM04}, and~\cite{Vishwanathy+etal:IEEE-ANTS09}. In this case, the average packet rate is nothing but the arrival rate of the Poisson process, which undergoes the upsurge triggered by an anomaly. Whether prior to or during an anomaly, the average packet rate typically measures in the thousands at the lowest, and therefore the parameter of the Poisson distribution is usually rather large. Hence, one can accurately approximate the Poisson distribution with a Gaussian one whose mean and variance are both equal to the average packet rate. This is precisely the $\mathcal{N}(\mu,a\mu)$-to-$\mathcal{N}(\theta,a\theta)$ model with $a=1$. However, in general, behavior of real traffic may deviate from the Poisson model, and parameter $a$ in the $\mathcal{N}(\mu,a\mu)$-to-$\mathcal{N}(\theta,a\theta)$ model can be used as an extra ``degree of freedom'' allowing one to take these deviations into account.  Thus, the $\mathcal{N}(\mu,a\mu)$-to-$\mathcal{N}(\theta,a\theta)$ model is a more general and better option.

The main objective of this work is to carry out a peer-to-peer comparative multiple-measure performance analysis of four detection procedures: the CUSUM procedure, the Shiryaev--Roberts (SR) procedure, and two of its derivatives -- the Shiryaev--Roberts--Pollak (SRP) procedure and the Shiryaev--Roberts--$r$ (SR--$r$) procedure. The performance measures we are interested in are Pollak's~\citeyearpar{Pollak:AS85} Supremum Average Detection Delay (ADD) and Shiryaev's~\citeyearpar{Shiryaev:TPA63} Stationary ADD -- each subject to a tolerable lower bound, $\gamma$, on the Average Run Length (ARL) to false alarm. Our intent is to test the asymptotic (as $\gamma\to\infty$) approximations for each performance index and each procedure of interest.

The remainder of the paper is organized as follows. In~\autoref{sec:opt-criteria+det-schemes} we formally state the problem, introduce the performance measures and the four procedures; for those with exact optimality properties we also remind these properties. In~\autoref{sec:asymptotics} we discuss the asymptotic properties of the detection procedures and present the corresponding asymptotic performance approximations. \autoref{s:perf-eval} describes the methodology of evaluating operating characteristics.  \autoref{sec:analysis} is devoted to the $\mathcal{N}(\mu,a\mu)$-to-$\mathcal{N}(\theta,a\theta)$ change-point scenario. Specifically, using the methodology of~\autoref{s:perf-eval}, we study the performance of each of the four procedures and examine the accuracy of the corresponding asymptotic approximations. \autoref{sec:cybersecurity} is intended to illustrate how change-point detection in general and the $\mathcal{N}(\mu,a\mu)$-to-$\mathcal{N}(\theta,a\theta)$ model in particular can be used in cybersecurity for rapid anomaly detection in computer networks. Lastly,~\autoref{sec:conclusion} draws conclusions.

\section{Optimality criteria and detection procedures}
\label{sec:opt-criteria+det-schemes}

As indicated in the introduction, the focus of this work is on two formulations of the change-point detection problem -- the minimax formulation and that related to multi-cyclic disorder detection in a stationary regime. The two, though both stem from the same assumption that the change-point is an unknown (non-random) number, posit their own optimality criterion.

Suppose there is an ``experimenter'' who is able to gather, one by one, a series of independent random observations, $\{X_n\}_{n\ge1}$. Statistically, the series is such that, for some time index $\nu$, which is referred to as the change-point, $X_1,X_2,\ldots, X_{\nu}$ each posses a known pdf $f(x)$, and $X_{\nu+1},X_{\nu+2},\ldots$ are each drawn from a population with a pdf $g(x)\not\equiv f(x)$, also known. The change-point is the unknown serial number of the last $f(x)$-distributed observation; therefore, if $\nu=\infty$, then the entire series $\{X_n\}_{n\ge1}$ is sampled from distribution with the pdf $f(x)$, and if $\nu=0$, then all observations are $g(x)$-distributed. The experimenter's objective is to decide that the change is in effect, raise an alarm and respond. The challenge is to make this decision ``as soon as possible'' past and ``no earlier'' than a certain prescribed limit prior to the true change-point.

Statistically, the problem is to sequentially differentiate between the hypotheses $\mathcal{H}_k\colon\nu=k\ge0$, i.e., that the change occurs at time moment $\nu=k$, $0\le k<\infty$, and $\mathcal{H}_{\tinyinfty}\colon\nu=\infty$, i.e., that the change never strikes. Once the $k$-th observation is made, the experimenter's decision options are either to accept $\mathcal{H}_k$, and thus declare that the change has occurred, or to reject $\mathcal{H}_k$ and continue observing data.

To test $\mathcal{H}_k$ against $\mathcal{H}_{\tinyinfty}$, one first constructs the corresponding likelihood ratio (LR). Let $\boldsymbol{X}_{1:n}=(X_1,X_2,\ldots,X_n)$ be the vector of the first $n\ge1$ observations; then the joint pdf-s of $\boldsymbol{X}_{1:n}$ under $\mathcal{H}_k$ and $\mathcal{H}_{\tinyinfty}$ are
\begin{align*}
p(\boldsymbol{X}_{1:n}|\mathcal{H}_{\tinyinfty})
&=
\prod_{j=1}^n f(X_j)
\quad\text{and}\quad
p(\boldsymbol{X}_{1:n}|\mathcal{H}_k)
=
\Biggl(\,\prod_{j=1}^{k} f(X_j)\Biggr)\times\Biggl(\,\prod_{j=k+1}^n g(X_j)\Biggr),
\end{align*}
where hereafter it is understood that $\prod_{j=k+1}^n=1$ if $k\ge n$; i.e., $p(\boldsymbol{X}_{1:n}|\mathcal{H}_{\tinyinfty})=p(\boldsymbol{X}_{1:n}|\mathcal{H}_k)$ if $k\ge n$. For the LR, $\LR_{k:n}=p(\boldsymbol{X}_{1:n}|\mathcal{H}_k)/p(\boldsymbol{X}_{1:n}|\mathcal{H}_{\tinyinfty})$, one then obtains
\begin{align}\label{eq:LR-def}
\LR_{k:n}
&=
\prod_{j=k+1}^n\LR_j,
\;\;\text{where}\;\;
\LR_n
=
\frac{g(X_n)}{f(X_n)},
\end{align}
and we note that $\LR_{k:n}=1$ if $k\ge n$. We also assume that $\Lambda_0=1$.

Next, to decide which of the hypotheses $\mathcal{H}_k$ or $\mathcal{H}_{\tinyinfty}$ is true, the sequence $\{\LR_{k:n}\}_{1\le k\le n}$ is turned into a detection statistic. To this end, one can either go with the maximum likelihood principle and use the detection statistic $W_n=\max_{1\le k\le n}\LR_{k:n}$ (maximal LR) or with the generalized Bayesian approach (limit of the Bayesian approach) and use the quasi-Bayesian detection statistic $R_n=\sum_{k=1}^n\LR_{k:n}$ (average LR with respect to an improper uniform prior distribution of the change-point). See~\cite{Polunchenko+Tartakovsky:MCAP11} for an overview of change-point detection approaches.

Once the detection statistic is chosen, it is supplied to an appropriate sequential detection procedure. Given the series $\{X_n\}_{n\ge1}$, a detection procedure is defined as a stopping time, $\T$, adapted to the filtration $\{\mathcal{F}_n\}_{n\ge0}$, where $\mathcal{F}_n=\sigma(\boldsymbol{X}_{1:n})$ is the sigma-algebra generated by the observations collected up to time instant $n\ge1$. The meaning of $\T$ is that after observing $X_1,X_2,\ldots,X_{\T}$ it is declared that the change is in effect. That may or may not be the case. If it is not, then $\T\le\nu$, and it is said that a false alarm has been sounded.

The above two detection statistics -- $\{W_n\}_{n\ge 1}$ and $\{R_n\}_{n\ge 1}$ -- give raise to a myriad of detection procedures. The first one we will be interested in is the CUmulative SUM (CUSUM) procedure. This is a maximum LR-based ``inspection scheme'' proposed by~\cite{Page:B54} to help solve issues arising in the area of industrial quality control. Currently, the CUSUM procedure is the {\it de facto} ``workhorse'' in a number of branches of engineering. Formally, we define the CUSUM procedure as the stopping time
\begin{align}\label{eq:T-CS-def}
\mathcal{C}_A
&=
\inf\{n\ge1\colon W_n\ge A\},
\end{align}
where $A>0$ is the detection threshold, and $W_n=\max_{1\le k\le n}\LR_{k:n}$ is  the CUSUM statistic, which can also be defined recursively as
\begin{align}\label{eq:Wn-CS-def}
W_n
&=
\max\{1,W_{n-1}\}\LR_n,
\;\;n\ge1
\;\;\text{with}\;\; W_0=1,
\end{align}
Hereafter in the definitions of stopping times  $\inf\{\varnothing\}=\infty$.

Another detection procedure we will consider is the the Shiryaev--Roberts (SR) procedure. In contrast to the CUSUM procedure, the SR procedure is based on the quasi-Bayesian argument. It is due to the independent work of~\cite{Shiryaev:SMD61,Shiryaev:TPA63} and \cite{Roberts:T66}. Specifically, Shiryaev introduced this procedure in continuous time in the context of detecting a change in the drift of a Brownian motion, and Roberts addressed the discrete-time problem of detecting a shift in the mean of a sequence of independent Gaussian random variables.  Formally, the SR procedure is defined as the stopping time
\begin{align}\label{eq:T-SR-def}
\mathcal{S}_A
&=
\inf\{n\ge1\colon R_n\ge A\},
\end{align}
where $A>0$ is the detection threshold, and $R_n=\sum_{k=1}^n\LR_{k:n}$ is the SR statistic, which can also be computed recursively as
\begin{align}\label{eq:Rn-SR-def}
R_n
&=
(1+R_{n-1})\LR_n,\;\; n\ge1\;\;\text{with}\;\; R_0=0,
\end{align}
Note that  the SR statistic starts from {\em zero}.

We will also be interested in two derivatives of the SR procedure: the Shiryaev--Roberts--Pollak (SRP) procedure and the Shiryaev--Roberts--$r$ (SR--$r$) procedure. The former is due to~\cite{Pollak:AS85} whose idea was to start the SR detection statistic, $\{R_n\}_{n\ge0}$, not from zero, but from a random point $R_0=R_0^Q$, where $R_0^Q$ is sampled from the quasi-stationary distribution $Q_A(x)=\lim_{n\to\infty}\Pr_{\tinyinfty}(R_n\le x|\mathcal{S}_A>n)$ of the SR detection statistic under the hypothesis $\mathcal{H}_{\tinyinfty}$. (Note that $\{R_n\}_{n\ge0}$ is a Markov Harris-recurrent process under $\mathcal{H}_{\tinyinfty}$.)  Formally, the SRP procedure is defined as the stopping time
\begin{align}\label{eq:T-SRP-def}
\mathcal{S}_A^{Q}
&=
\inf\{n\ge1\colon R_n^{Q}\ge A\},
\end{align}
where $A>0$ is the detection threshold, and $R_n^Q$ is the SRP detection statistic given by the recursion
\begin{align}\label{eq:R-SRP-def}
R_{n}^{Q}
&=
(1+R_{n-1}^{Q})\LR_n,
\;\; n\ge1
\;\;\text{with}\;\;
R_0^{Q}\propto Q_A(x) .
\end{align}

The SR--$r$ procedure was proposed by~\cite{Moustakides+etal:SS11} who regard starting off the original SR detection statistic, $\{R_n\}_{n\ge0}$, from a fixed (but specially designed) $R_0=r$, $0\le r<A$, and defining the stopping time with this new deterministic initialization as
\begin{align}\label{eq:T-SR-r-def}
\mathcal{S}_A^r
&=
\inf\{n\ge1\colon R_n^r\ge A\},
\end{align}
where $A>0$ is the detection threshold and the SR--$r$ detection statistic $R_n^r$ is given by the recursion
\begin{align}\label{eq:R-SR-r-def}
R_n^r
&=
(1+R_{n-1}^r)\LR_n,
\;\; n\ge1~~\text{with}\;\;
R_0^r=r\ge0,
\end{align}
Note that for $r=0$ this is nothing but the conventional SR procedure.

We now proceed with reviewing the optimality criteria whereby one decides which procedure to use.

We first set down some additional notation. Let $\Pr_{\nu}(\cdot)$ be the probability measure generated by the observations $\{X_n\}_{n\ge1}$ when the change-point is $\nu\ge0$, and $\EV_{\nu}[\,\cdot\,]$ be the corresponding expectation. Likewise, let $\Pr_{\tinyinfty}(\cdot)$ and $\EV_{\tinyinfty}[\,\cdot\,]$ denote the same under the no-change scenario, i.e., when $\nu=\infty$.

Consider first  the minimax formulation proposed by~\cite{Lorden:AMS71} where the risk of raising a false alarm is measured by the ARL to false alarm $\ARL(\T)=\EV_{\tinyinfty}[\T]$ and the delay to detection by the ``worst-worst-case'' ADD $\mathcal{J}_{\mathrm{L}}(\T)=\sup_{0\le\nu<\infty}\{\esssup\EV_{\nu}[(\T-\nu)^+|\mathcal{F}_{\nu}]\}$, where hereafter $x^+=\max\{0,x\}$. Let $\Delta(\gamma)=\{\T\colon\ARL(\T)\ge\gamma\}$ be the class of detection procedures (stopping times) for which the ARL to false alarm does not fall below a given (desired and {\it a~priori} set) level $\gamma>1$.  Lorden's version of the minimax optimization problem is to find $\T_{\mathrm{opt}}\in\Delta(\gamma)$ such that $\mathcal{J}_{\mathrm{L}}(\T_{\mathrm{opt}})=\inf_{\T\in\Delta(\gamma)}\mathcal{J}_{\mathrm{L}}(\T)$ for every $\gamma>1$. For the iid model, this problem was solved by~\cite{Moustakides:AS86}, who showed that the solution is the CUSUM procedure. Specifically, if we select the detection threshold $A=A_\gamma$ from the solution of the equation $\ARL(\mathcal{C}_{A_\gamma})=\gamma$, then $\mathcal{J}_{\mathrm{L}}(\mathcal{C}_{A_\gamma})=\inf_{\T\in\Delta(\gamma)}\mathcal{J}_{\mathrm{L}}(\T)$ for every $\gamma>1$.

Though CUSUM's strict $\mathcal{J}_{\mathrm{L}}(\T)$-optimality is a strong result, ideal for engineering purposes would be to have a procedure that minimizes the average (conditional) detection delay, $\ADD_{\nu}(\T)=\EV_{\nu}[\T-\nu|\T>\nu]$, for all $\nu\ge0$ simultaneously. As no such uniformly optimal procedure is possible,~\cite{Pollak:AS85}
suggested to revise Lorden's version of minimax optimality by replacing $\mathcal{J}_{\mathrm{L}}(\T)$ with $\mathcal{J}_{\mathrm{P}}(\T)=
\sup_{0\le\nu<\infty}\ADD_{\nu}(\T)$, i.e., with the worst average (conditional) detection delay, and seek $\T_{\mathrm{opt}}\in\Delta(\gamma)$ such that $\mathcal{J}_{\mathrm{P}}(\T_{\mathrm{opt}})=\inf_{\T\in\Delta(\gamma)}\mathcal{J}_{\mathrm{P}}(\T)$ for every $\gamma>1$. We believe that this problem has higher applied potential than Lorden's. Contrary to the latter, an exact solution to Pollak's minimax problem is still an open question. See, e.g.,~\cite{Pollak:AS85},~\cite{Polunchenko+Tartakovsky:AS10}, \cite{Moustakides+etal:SS11}, and \cite{Tartakovsky+etal:TPA11} for a related study.

Note that Lorden's and Pollak's versions of the minimax formulation both assume that the detection procedure is applied only once; the result is either a false alarm, or a correct (though delayed) detection, and no data sampling is done past the stopping point. This is known as the single-run paradigm. Yet another formulation emerges if one considers applying the same procedure in cycles, e.g., starting anew after every false alarm. This is the multi-cyclic formulation.

Specifically, the idea is to assume that in exchange for the assurance that the change will be detected with maximal speed, the experimenter agrees to go through a ``storm'' of false alarms along the way. The false alarms are ensued from repeatedly applying the same detection rule, starting from scratch after each false alarm. Put otherwise, suppose the change-point, $\nu$,  is substantially larger than the desired level of the ARL to false alarm, $\gamma>1$. That is, the change occurs in a distant future and it is preceded by a stationary flow of false alarms; the ARL to false alarm in this context is the mean time  between consecutive false alarms. As argued by~\cite{Pollak+Tartakovsky:SS09}, this comes in handy in many surveillance applications, in particular in the area of cybersecurity which will be addressed in~\autoref{sec:cybersecurity}.

Formally, let $T_1,T_2,\ldots$ denote sequential independent applications of the same stopping time $\T$, and let $\mathcal{T}_{(j)}=T_{(1)}+T_{(2)}+\cdots+T_{(j)}$ be the time of the $j$-th alarm, $j\ge1$. Let $I_\nu=\min\{j\ge1\colon\mathcal{T}_{(j)}>\nu\}$ so that $\mathcal{T}_{(I_\nu)}$ is the point of detection of the true change, which occurs at time instant $\nu$ after $I_\nu-1$ false alarms have been raised. Consider $\mathcal{J}_{\mathrm{ST}}(\T)=\lim_{\nu\to\infty}\EV_\nu[\mathcal{T}_{(I_\nu)}-\nu]$, i.e., the limiting value of the ADD that we will refer to as the {\em stationary ADD} (STADD); then the multi-cyclic optimization problem consists in finding $\T_{\mathrm{opt}}\in\Delta(\gamma)$ such that $\mathcal{J}_{\mathrm{ST}}(\T_{\mathrm{opt}})=
\inf_{\T\in\Delta(\gamma)}\mathcal{J}_{\mathrm{ST}}(\T)$ for every $\gamma>1$. For the continuos-time Brownian motion model, this formulation was first proposed by  \cite{Shiryaev:SMD61,Shiryaev:TPA63} who showed that the SR procedure is strictly optimal. For the discrete-time iid model, optimality of the SR procedure in this setting was solved recently by~\cite{Pollak+Tartakovsky:SS09}. Specifically, introduce the integral average detection delay
\begin{align}\label{eq:IADD-def}
\IADD(\T)
&=
\sum_{\nu=0}^{\infty}\EV_{\nu}[(\T-\nu)^+],
\end{align}
and the relative integral average detection delay $\mathcal{J}_{\mathrm{GB}}(T)=\IADD(\T)/\ARL(\T)$. If the detection threshold, $A$, is set to $A_\gamma$, the solution of the equation $\ARL(\mathcal{S}_{A_\gamma})=\gamma$,  then $\mathcal{J}_{\mathrm{GB}}(\mathcal{S}_{A_\gamma})=\inf_{\T\in\Delta(\gamma)}\mathcal{J}_{\mathrm{GB}}(\T)$ for every $\gamma>1$. Also $\mathcal{J}_{\mathrm{GB}}(\T)\equiv\mathcal{J}_{\mathrm{ST}}(\T)$ for any stopping time $\T$, so that $\mathcal{J}_{\mathrm{ST}}(\mathcal{S}_{A_\gamma})=\inf_{\T\in\Delta(\gamma)}\mathcal{J}_{\mathrm{ST}}(\T)$ for every $\gamma>1$.

We conclude this section reiterating that our goal is to test the above four procedures -- CUSUM, SR, SRP and SR--$r$ -- against each other with respect to two measures of detection delay -- $\mathcal{J}_{\mathrm{P}}(\T)$ and $\mathcal{J}_{\mathrm{ST}}(\T)$. We will also analyze the accuracy of the asymptotic approximations which are introduced in the next section.

\section{Asymptotic optimality and performance approximations}
\label{sec:asymptotics}

To remind, we focus on Pollak's minimax formulation of minimizing the maximal ADD and on the multi-cyclic formulation of minimizing the stationary ADD. We also indicated that the latter is a solved problem (and the solution is the SR procedure), while the former is still an open question. The usual way around is to consider the asymptotic case $\gamma\to\infty$. The hope is to design a procedure, $\T_{\mathrm{opt}}^*\in\Delta(\gamma)$, such that $\mathcal{J}_{\mathrm{P}}(\T_{\mathrm{opt}}^*)$ and the (unknown) optimum, $\inf_{\T\in\Delta(\gamma)}\mathcal{J}_{\mathrm{P}}(\T)$, will be in some sense ``close'' to each other in the limit, as $\gamma\to\infty$. To this end, three different types of $\mathcal{J}_{\mathrm{P}}(\T_{\mathrm{opt}}^*)$-to-$\inf_{\T\in\Delta(\gamma)}\mathcal{J}_{\mathrm{P}}(\T)$ convergence are generally distinguished:
\begin{inparaenum}[\itshape a\upshape)]
    \item A procedure $\T_{\mathrm{opt}}^*\in\Delta(\gamma)$ is said to be first-order asymptotically $\mathcal{J}_{\mathrm{P}}(\T)$-optimal in the class $\Delta(\gamma)$ if $\mathcal{J}_{\mathrm{P}}(\T_{\mathrm{opt}}^*)=[\inf_{\T\in\Delta(\gamma)}\mathcal{J}_{\mathrm{P}}(\T)][1+o(1)]$, as $\gamma\to\infty$, where hereafter $o(1)\to0$, as $\gamma\to\infty$,
    \item second-order -- if
$\mathcal{J}_{\mathrm{P}}(\T_{\mathrm{opt}}^*)=\inf_{\T\in\Delta(\gamma)}\mathcal{J}_{\mathrm{P}}(\T)+O(1)$, as $\gamma\to\infty$, where $O(1)$ stays bounded, as $\gamma\to\infty$, and
    \item third-order -- if $\mathcal{J}_{\mathrm{P}}(\T_{\mathrm{opt}}^*)=\inf_{\T\in\Delta(\gamma)}\mathcal{J}_{\mathrm{P}}(\T)+o(1)$, as $\gamma\to\infty$.
\end{inparaenum}
Note that $\inf_{\T\in\Delta(\gamma)}\mathcal{J}(\T)\to\infty$ as $\gamma\to\infty$.

We now review our four procedures' individual asymptotic optimality properties and provide the corresponding asymptotic approximations for $\ARL(\T)=\EV_{\tinyinfty}[\T]$, $\mathcal{J}_{\mathrm{P}}(\T)=\sup_{0\le\nu<\infty}\ADD_{\nu}(\T)$ (where $\ADD_{\nu}(\T)=\EV_{\nu}[\T-\nu|\T>\nu]$),   $\ADD_{\tinyinfty}(\T)=\lim_{\nu\to\infty}\ADD_{\nu}(\T)$, and $\mathcal{J}_{\mathrm{ST}}(T)$.  Let $Z_i=\log\LR_i$ and let $I_f=-\EV_{\tinyinfty}[Z_1]$ and $I_g=\EV_0[Z_1]$ denote the Kullback--Leibler information numbers. Henceforth, it will be assumed that $Z_1$ is $\Pr_{\tinyinfty}$- and $\Pr_0$-nonarithmetic and that $0<I_f<\infty$ and $0<I_g<\infty$. Throughout the rest of the paper we will also assume the second moment conditions $\EV_\infty |Z_1|^2<\infty$ and $\EV_0 |Z_1|^2<\infty$. Further, introduce the random walk $\{S_n\}_{n\ge0}$, where $S_n=\sum_{j=1}^n Z_i$, $n\ge1$, with $S_0=0$. Let $\tilde{V}_{\tinyinfty}=\sum_{j=1}^{\infty} e^{-S_j}$, and define $\tilde{Q}(x)=\Pr_0(\tilde{V}_{\tinyinfty}\le x)$. Next, for $a\ge 0$, introduce the one-sided stopping time $\tau_a=\inf\{n\ge1\colon S_n\ge a\}$ and let $\kappa_a=S_{\tau_a}-a$ denote the overshoot (i.e., the excess of $S_n$ over the level $a$ at stopping). Define $\zeta= \lim_{a\to\infty}\EV_0[e^{-\kappa_a}]$ and $\varkappa=\lim_{a\to\infty}\EV_0[\kappa_a]$. It can be shown that
\begin{align}\label{eq:zeta+kappa}
\zeta
&=
\frac{1}{I_g}\exp\left\{-\sum_{k=1}^\infty\frac{1}{k}\bigl[\Pr_{\tinyinfty}(S_k>0)+\Pr_0(S_k\le 0)\bigr]\right\}
\;\;\text{and}\;\;
\varkappa
=
\frac{\EV_0[Z_1^2]}{2\EV_0[Z_1]}-\sum_{k=1}^\infty\frac{1}{k}\EV_0[S_k^-],
\end{align}
where hereafter $x^-=-\min(0,x)$; cf., e.g.,~\citet[Chapters~2~\&~3]{Woodroofe:Book82} and~\citet[Chapter~VIII]{Siegmund:Book85}.

We are now poised to discuss  optimality properties and approximations for the operating characteristics of the four detection procedures of interest. We begin with the CUSUM procedure. First, from the fact that $\mathcal{J}_{\mathrm{L}}(\mathcal{C}_A)=\mathcal{J}_{\mathrm{P}}(\mathcal{C}_A)= \ADD_0(\mathcal{C}_A)$ and the work of~\cite{Lorden:AMS71} one can deduce that the CUSUM procedure is first-order asymptotically $\mathcal{J}_{\mathrm{P}}(\T)$-optimal. However, in fact it is second-order $\mathcal{J}_{\mathrm{P}}(\T)$-minimax which can be established as follows.  Since CUSUM is strictly optimal in the sense of minimizing the essential supremum ADD $\mathcal{J}_{\rm L}(T)$, which is equal to $\ADD_0$ for both CUSUM and SR, it is clear that  $\ADD_0(\mathcal{C}_A) < \ADD_0(\mathcal{S}_A)$, where the thresholds are different for  CUSUM and SR to attain the same ARL to false alarm $\gamma$. By  \cite{Tartakovsky+etal:TPA11}, the SR procedure is second-order asymptotically minimax with respect to $\mathcal{J}_{\mathrm{P}}(T)$, so that CUSUM is also second-order asymptotically $\mathcal{J}_{\mathrm{P}}(\T)$-optimal.  Now, let $A=A_\gamma$, where $A_\gamma$ is the solution of the equation $\ARL(\mathcal{C}_{A_\gamma})=\gamma$. Then
\begin{equation} \label{eq:SADDCS}
\mathcal{J}_{\mathrm{P}}(\mathcal{C}_{A_\gamma})=
\frac{1}{I_g}(\log A_\gamma+\varkappa+\beta_0)+o(1),\;\;\text{as}\;\;\gamma\to\infty,
\end{equation}
where $\beta_0=\EV_0[\min_{n\ge0} S_n]$. We iterate that $\mathcal{J}_{\mathrm{P}}(\mathcal{C}_{A})=\EV_0[\mathcal{C}_{A}]$. This asymptotic expansion was first obtained by~\cite{Dragalin:PSIM94} for the single-parameter exponential family, but it holds in a general case too as long as the second moment condition $\EV_0[Z_1^2]<\infty$ is satisfied and $Z_1$ is $\Pr_0$-non-arithmetic. See~\cite{Tartakovsky:IEEE-CDC05}, who also shows that
\begin{equation} \label{eq:ADDinftyCS}
\ADD_{\tinyinfty}(\mathcal{C}_{A_\gamma}) =
\frac{1}{I_g}(\log A_\gamma+\varkappa-\beta_{\tinyinfty})+o(1),\;\;\text{as}\;\;\gamma\to\infty,
\end{equation}
where $\beta_{\tinyinfty}=\lim_{n\to\infty}\EV_{\tinyinfty}[S_n-\min_{0\le k\le n}S_k]$; constants $\beta_0$ and $\beta_{\tinyinfty}$ can be computed numerically (e.g., by Monte Carlo simulations).  An accurate approximation for the ARL to false alarm is as follows:
\begin{equation}\label{eq:ARLcusum}
\ARL(\mathcal{C}_{A}) \approx \frac{A}{I_g\zeta^2}-\frac{\log A}{I_f}-\frac{1}{I_g\zeta} .
\end{equation}

We now proceed to the SR procedure. First, recall that, by~\cite{Pollak+Tartakovsky:SS09}, this procedure is a strictly optimal multi-cyclic procedure in the sense of minimizing the stationary ADD.  The following lower bound for the minimal SADD is fundamental for the analysis carried out in~\autoref{sec:analysis}: for every $r\ge0$,
\begin{align}\label{eq:LowerBound-def}
\inf_{\T\in\Delta(\gamma)}\mathcal{J}_{\mathrm{P}}(\T)\ge \mathcal{J}_{\mathrm{LB}}(\mathcal{S}_{A_\gamma}^r)
&=
\cfrac{r\ADD_0[\mathcal{S}_{A_\gamma}^r]+\IADD(\mathcal{S}_{A_\gamma}^r)}{r+\ARL(\mathcal{S}_{A_\gamma}^r)},
\end{align}
where $A_\gamma$ is the solution of the equation $\ARL(\mathcal{S}_{A_\gamma}^r)=\gamma$, $\gamma>1$
(see~\citealt{Moustakides+etal:SS11} and~\citealt{Polunchenko+Tartakovsky:AS10}). Taking $r=0$ in~\eqref{eq:LowerBound-def}, we obtain $\inf_{\T\in\Delta(\gamma)}\mathcal{J}_{\mathrm{P}}(\T)\ge  \mathcal{J}_{\mathrm{ST}}(\mathcal{S}_{A_\gamma})$. Let $R_{\tinyinfty}$ be a random variable that has the $\Pr_{\tinyinfty}$-limiting (stationary) distribution of $R_n$ as $n\to\infty$, i.e., $Q_{\mathrm{ST}}(x)=\lim_{n\to\infty}\Pr_{\tinyinfty}(R_n\le x)=\Pr_{\tinyinfty}(R_{\tinyinfty}\le x)$. Let $V = \sum_{k=1}^\infty e^{-S_k}$ and $\tilde{Q}(x) = \Pr_0 (V \le x)$. Let
\begin{align}\label{eq:C-inf-def}
C_{\tinyinfty}
&=
\EV[\log(1+R_{\tinyinfty}+\tilde{V})]
=
\int_0^\infty\int_0^\infty\log(1+x+y)\,dQ_{\mathrm{ST}}(x)\,d\tilde{Q}(y).
\end{align}
By~\cite{Tartakovsky+etal:TPA11},
\[
\mathcal{J}_{\mathrm{ST}}(\mathcal{S}_{A_\gamma}) = \frac{1}{I_g}(\log A_\gamma+\varkappa -C_\infty)+o(1),\;\;\text{as}\;\;\gamma\to\infty,
\]
so that
\begin{equation}\label{LBasympt}
\inf_{\T\in\Delta(\gamma)}\mathcal{J}_{\mathrm{P}}(\T)\ge \frac{1}{I_g}(\log A_\gamma+\varkappa -C_\infty)+o(1),\;\;\text{as}\;\;\gamma\to\infty .
\end{equation}
On the other hand, a straightforward argument shows that, for the SR procedure, if $A=A_\gamma$ is the solution of the equation $\ARL(\mathcal{S}_{A_\gamma})=\gamma$,  then $\ARL(\mathcal{S}_{A_\gamma})= \gamma \approx A_\gamma/\zeta$ and
\begin{align} \label{eq:SADDSR}
\mathcal{J}_{\mathrm{P}}(\mathcal{S}_{A_\gamma}) = \EV_0[\mathcal{S}_{A_\gamma}] =
\frac{1}{I_g}(\log A_\gamma+\varkappa-C_0)+o(1),\;\;\text{as}\;\;\gamma\to\infty,
\end{align}
where
\begin{align*}
C_0
&=
\EV[\log(1+\tilde{V}_{\tinyinfty})] = \int_0^\infty\log(1+x)\,d\tilde{Q}(x)
\end{align*}
(cf.~\cite{Tartakovsky+etal:TPA11}). Since $C_0< C_\infty$, it follows that the SR procedure is second-order minimax, the fact that we have already mentioned above when considering CUSUM. The difference is $(C_\infty-C_0)/I_g$, which can be quite large when detecting small changes.

For an arbitrary scenario, the constant $C_0$ and distribution $\tilde{Q}(x)$ are amenable to numerical treatment. For cases where both can be computed exactly in a closed form see~\cite{Tartakovsky+etal:TPA11} and~\cite{Polunchenko+Tartakovsky:MCAP11}. The approximation $\ARL(\mathcal{C}_A)\approx A/\zeta$ is known to be very accurate and is therefore recommended for practical use; it can be derived from the  the fact that $\{R_n-n\}_{n\ge0}$ is a zero-mean $\Pr_{\tinyinfty}$-martingale.

We now consider the SRP procedure. From decision theory (see, e.g.,~\citealp[Theorem~2.11.3]{Ferguson:Book67}) it is known that the minimax procedure should be an equalizer, i.e., $\ADD_{\nu}(\T)$ should be the same for all $\nu\ge0$. To make the SR procedure an equalizer, \cite{Pollak:AS85}  suggested to start the SR detection statistic, $\{R_n\}_{n\ge0}$, from a random point distributed according to the quasi-stationary distribution.  However, Pollak was able to demonstrate that the SRP procedure is only asymptotically third-order $\mathcal{J}_{\mathrm{P}}(\T)$-optimal. More specifically, let $\EV_0[Z_1^+]<\infty$, and suppose that the detection threshold, $A$, is set to $A_{\gamma}$, the solution of the equation $\ARL(\mathcal{S}_{A_\gamma}^Q)=\gamma$. Then $\mathcal{J}_{\mathrm{P}}(\mathcal{S}_{A_\gamma}^Q)=
\inf_{\T\in\Delta(\gamma)}\mathcal{J}_{\mathrm{P}}(\T)+o(1)$, as
$\gamma\to\infty$. Recently,~\cite{Tartakovsky+etal:TPA11} obtained an asymptotic approximation for $\mathcal{J}_{\mathrm{P}}(\mathcal{S}_A^Q)$ under the second moment condition $\EV_0[Z_1^2]<\infty$:
\begin{align}\label{eq:SADDSRP}
\mathcal{J}_{\mathrm{P}}(\mathcal{S}_{A_\gamma}^Q)
=
\frac{1}{I_g}(\log A_\gamma+\varkappa-C_{\tinyinfty})+o(1),\;\;\text{as}\;\;\gamma\to\infty,
\end{align}
and we note that $\mathcal{J}_{\mathrm{P}}(\mathcal{S}_{A}^Q)=\ADD_{\nu}(\mathcal{S}_{A}^Q)$ for all $\nu\ge0$, and $\mathcal{J}_{\mathrm{P}}(\mathcal{S}_{A}^Q)=\mathcal{J}_{\mathrm{ST}}(\mathcal{S}_{A}^Q)$. To approximate $\ARL(\mathcal{S}_{A}^Q)$,~\cite{Tartakovsky+etal:TPA11} suggest to use the formula $\ARL(\mathcal{S}_{A}^Q)\approx A/\zeta-\mu_Q$, where $\mu_Q= \int_0^A x\,dQ_A(x)$ is  the mean of the quasi-stationary distribution.

It is left to consider the SR--$r$ procedure. Since this procedure is sensitive to the choice of the starting point, we first discuss the question of how to design this point. To this end, fundamental is the lower bound for the minimal SADD given in \eqref{eq:LowerBound-def}.
From~\eqref{eq:LowerBound-def} one can deduce that if the starting point $R_0^r=r^*$ is chosen so that the SR--$r$ procedure is an equalizer (i.e., $\ADD_{\nu}(\mathcal{S}_A^{r^*})$ is the same for all $\nu\ge0$), then it will be exactly $\mathcal{J}_{\mathrm{P}}(\T)$-optimal. Indeed, if the SR--$r$ procedure is an equalizer, then $\mathcal{J}_{\mathrm{LB}}(\mathcal{S}_{A}^{r*})=\EV_0[\mathcal{S}_A^{r^*}]=\mathcal{J}_{\mathrm{P}}(\mathcal{S}_A^{r^*})$. This argument was exploited by~\cite{Polunchenko+Tartakovsky:AS10}, who found  change-point scenarios where $\ARL(\mathcal{S}_A^r)$,  $\ADD_{\nu}(\mathcal{S}_A^r)$, $\nu\ge0$, and $r^*$ can be obtained in a closed form and the SR$-r$ procedure with $r=r^*$ is strictly optimal. In general,~\cite{Moustakides+etal:SS11} suggest to set the starting point to the solution of the following constraint minimization problem:
\begin{equation}\label{eq:rstar}
r^*= \arginf_{0\le r<\infty} \left\{\mathcal{J}_{\mathrm{P}}(\mathcal{S}_A^r)-\mathcal{J}_{\mathrm{LB}}(\mathcal{S}_A^r)\right\}
\;\;\text{subject to}\;\; \ARL(\mathcal{S}_A^r)=\gamma .
\end{equation}
We discuss this approach in detail in~\autoref{sec:analysis} where we also discuss certain properties of $r^*$, as well as  another approach that allows one to obtain a nearly optimal result for the low false alarm rate.

As shown by~\cite{Tartakovsky+etal:TPA11}, the SR--$r$ procedure also enjoys the third-order asymptotic optimality property.  Specifically, assume that the  detection threshold $A=A_\gamma$ is the solution of the equation $\ARL(\mathcal{S}_{A_\gamma}^r)=\gamma$ and the head start $r$ (either fixed or growing at the rate of $o(\gamma)$) is selected in such a way that $\ARL(\mathcal{S}_{A_\gamma}^r) \approx A_\gamma/\zeta$ and that the supremum ADD $\mathcal{J}_{\mathrm{P}}(\mathcal{S}_{A_\gamma}^r)$ is attained at infinity, i.e., $\mathcal{J}_{\mathrm{P}}(\mathcal{S}_{A_\gamma}^r)=\ADD_{\tinyinfty}(\mathcal{S}_{A_\gamma}^r)$. Then
\begin{equation}\label{eq:SADDSRr}
\mathcal{J}_{\mathrm{P}}(\mathcal{S}_{A_\gamma}^r)= \ADD_{\tinyinfty}(\mathcal{S}_{A_\gamma}^r)=
\frac{1}{I_g}(\log A_\gamma+\varkappa-C_{\tinyinfty})+o(1), \;\;\text{as}\;\;\gamma\to\infty.
\end{equation}
Comparing with the lower bound \eqref{LBasympt} we see that in this case the SR$-r$ procedure is nearly optimal (within the negligible term $o(1)$).  For practical purposes, to approximate $\ARL(\mathcal{S}_{A}^r)$ one may use the formula $\ARL(\mathcal{S}_{A}^r)\approx A/\zeta-r$, which is rather accurate, and can be easily obtained from the fact that $\{R_n^r-n-r\}_{n\ge0}$ is a zero-mean $\Pr_{\tinyinfty}$-martingale.

\section{Numerical techniques for evaluation of operating characteristics}
\label{s:perf-eval}

In this section we present integral equations for operating characteristics of change detection procedures as well as outline numerical techniques for solving these equations that can be effectively used for performance evaluation. Further details can be found in~\cite{Moustakides+etal:SS11}.

Consider a generic detection procedure described by the stopping time
\begin{align}\label{eq:generic-T-def}
\mathcal{T}_A^s
&=
\inf\{n\ge1\colon V_n^s\ge A\}, \quad A >0.
\end{align}
where $\{V_n^s\}_{n\ge0}$ is a generic Markov detection statistic that follows the recursion
\begin{align}\label{eq:generic-V-def}
V_{n}^s
&=
\xi(V_{n-1}^s)\LR_{n},
\;\;n\ge1
\;\;\text{with}\;\; V_0^s=s\ge0.
\end{align}
Here $\xi(x)$ is a non-negative function and $s$ is a fixed parameter, referred to as the starting point or the head start.

The above generic stopping time describes a rather broad class of detection procedures; in particular, all considered procedures belong to this class. Indeed, for the CUSUM procedure $\xi(x)=\max\{1,x\}$ and for the SR-type procedures $\xi(x)=1+x$. This universality of $\mathcal{T}_A^s$ enables one to evaluate practically any Operating Characteristic (OC) of any procedure that is a special case of $\mathcal{T}_A^s$, merely by picking the right $\xi(x)$. We now provide a set of exact integral equations on all major OC-s for $\mathcal{T}_A^s$.

Recall that $\LR_n=g(X_n)/f(X_n)$, $n\ge1$, and assume that $\LR_1$ is continuous. For $d=\{0,\infty\}$, let $P_d^{\LR}(t)=\Pr_d(\LR_1\le t)$ denote the cdf of the LR. Let
\begin{align*}
\mathcal{K}_d(x,y)
&=
\frac{\partial}{\partial y}\Pr_d(V_{n+1}^s\le y|V_n^s=x)
=
\frac{\partial}{\partial y}P_d^{\LR}\left(\frac{y}{\xi(x)}\right),
\; d=\{0,\infty\}
\end{align*}
denote the transition probability density kernel for the Markov process $\{V_n^s\}_{n\ge1}$. Note that $dP_0^{\LR}(t)=t\,dP_{\tinyinfty}^{\LR}(t)$, and therefore, $\mathcal{K}_0(x,y)/\mathcal{K}_{\tinyinfty}(x,y)=y/\xi(x)$, which can be used as a ``shortcut'' in deriving the formula for $\mathcal{K}_0(x,y)$ from that for $\mathcal{K}_{\tinyinfty}(x,y)$ or vice versa.

Suppose that $V_0^x=x$ is fixed. Let $\ell(x)=\EV_{\tinyinfty}[\mathcal{T}_A^x]$ and $\delta_0(x)=\EV_0[\mathcal{T}_A^x]$. It can be shown that
$\ell(x)$ and $\delta_0(x)$ satisfy the {\em renewal equations}
\begin{align}\label{eq:ARL+ADD0-int-eqn}
\ell(x)
&=
1+\int_0^A\mathcal{K}_{\tinyinfty}(x,y)\,\ell(y)\,dy
\;\;\text{and}\;\;
\delta_0(x)
=
1+\int_0^A\mathcal{K}_0(x,y)\,\delta_0(y)\,dy,
\end{align}
respectively (cf.~\cite{Moustakides+etal:SS11}). Next, consider $\ADD_{\nu}(\mathcal{T}_A^x)=\EV_\nu[\mathcal{T}_A^x-\nu|\mathcal{T}_A^x>\nu]$ for an arbitrary fixed $\nu\ge1$; note that $\ADD_0(\mathcal{T}_A^x)\equiv\delta_0(x)$. To evaluate $\ADD_{\nu}(\mathcal{T}_A^x)$,~\cite{Moustakides+etal:SS11} first argue that $\Pr_{\nu}(\mathcal{T}_A^x>\nu)=\Pr_{\tinyinfty}(\mathcal{T}_A^x>\nu)$, and consequently, $\ADD_{\nu}(\mathcal{T}_A^x)=\EV_{\nu}[(\mathcal{T}_A^x-\nu)^+]\,/\,\Pr_{\tinyinfty}(\mathcal{T}_A^x>\nu)$, so that we turn attention to $\delta_{\nu}(x)=\EV_{\nu}[(\mathcal{T}_A^x-\nu)^+]$ and $\rho_{\nu}(x)=\Pr_{\tinyinfty}(\mathcal{T}_A^x>\nu)$. It is direct to see that
\begin{align*}
\delta_{\nu}(x)
&=
\int_0^A\mathcal{K}_{\tinyinfty}(x,y)\,\delta_{\nu-1}(y)\,dy
\;\;\text{and}\;\;
\rho_{\nu}(x)=\int_0^A\mathcal{K}_{\tinyinfty}(x,y)\,\rho_{\nu-1}(y)\,dy, \quad \nu \ge 1,
\end{align*}
where  $\delta_0(x)$ is as in~\eqref{eq:ARL+ADD0-int-eqn} and $\rho_0(x)\equiv1$, since $\Pr_{\tinyinfty}(\mathcal{T}_A^x>0)\equiv1$; cf.~\cite{Moustakides+etal:SS11}. As soon as $\delta_{\nu}(x)$ and $\rho_{\nu}(x)$ are found, by the above argument $\ADD_{\nu}(\mathcal{T}_A^x)$ can be evaluated as the ratio $\delta_{\nu}(x)/\rho_{\nu}(x)$. Furthermore, using $\ADD_{\nu}(\mathcal{T}_A^x)$'s computed for sufficiently many successive $\nu$'s beginning from $\nu=0$ and higher, one can also evaluate $\mathcal{J}_{\mathrm{P}}(\mathcal{T}_A^x)=\sup_{0\le\nu<\infty}\ADD_{\nu}(\mathcal{T}_A^x)$, since $\ADD_{\tinyinfty}(\mathcal{T}_A^x)=\lim_{\nu\to\infty}\ADD_{\nu}(\mathcal{T}_A^x)$ is independent of $V_0^x=x$.

We now proceed to the stationary average detection delay $\mathcal{J}_{\mathrm{ST}}(\T)$. It follows from ~\cite{Pollak+Tartakovsky:SS09} that the STADD is equal to the relative integral ADD:
\begin{align*}
\mathcal{J}_{\mathrm{ST}}(\mathcal{T}_A^x)
&=
\left(\,\sum_{\nu=0}^{\infty}\EV_{\nu}[(\mathcal{T}_A^x-\nu)^+]\right)\bigm/\EV_{\tinyinfty}[\mathcal{T}_A^x],
\end{align*}
and if we let $\psi(x)=\sum_{\nu=0}^{\infty}\EV_{\nu}[(\mathcal{T}_A^x-\nu)^+]=\sum_{\nu=0}^{\infty}\delta_{\nu}(x)$, then $\mathcal{J}_{\mathrm{ST}}(\mathcal{T}_A^x)=\psi(x)/\ell(x)$. It can be shown that $\psi(x)$ satisfies
\begin{align}\label{eq:IADD-int-eqn}
\psi(x)
&=
\delta_0(x)+\int_0^A\mathcal{K}_{\tinyinfty}(x,y)\,\psi(y)\,dy,
\end{align}
where $\delta_0(x)$ is as in~\eqref{eq:ARL+ADD0-int-eqn}.

The IADD is also involved in the lower bound $\mathcal{J}_{\mathrm{LB}}(\T)$ which is defined in~\eqref{eq:LowerBound-def} and is specific to the SR--$r$ procedure, i.e., for the case when $\xi(x)=1+x$. Assuming this choice of $\xi(x)$, observe that $\mathcal{J}_{\mathrm{LB}}(\mathcal{S}_A^r)=[r\delta_0(r)+\psi(r)]/[r+\ell(r)]$, which can be evaluated once $\ell(x)$, $\delta_0(x)$, and $\psi(x)$ are found.

Consider now randomizing the starting point $V_0^x=x$ in a fashion analogous that behind the SRP stopping time. Let $Q_A(x)=\lim_{n\to\infty}\Pr_{\tinyinfty}(V_n^s\le x|\mathcal{T}_A^s>n)$ be quasi-stationary distribution; note that this distribution exists, as guaranteed by~\cite[Theorem~III.10.1]{Harris:Book63}. Let $\mathcal{T}_A^Q=\inf\{n\ge1\colon V_n^Q\ge A\}$, where  $\{V_n^Q\}_{n\ge0}$ is a randomized generic detection statistic such that $V_{n}^{Q}=\xi(V_{n-1}^{Q})\LR_n$ for $n\ge1$ with $V_0^Q\propto Q_A$. Note that $\mathcal{T}_A^Q$ is the SRP procedure when $\xi(x)=1+x$.

For $\mathcal{T}_A^Q$, any OC is dependent upon the quasi-stationary distribution. We therefore first state the equation that determines the quasi-stationary pdf $q_A(x)=dQ_A(x)/dx$:
\begin{align}\label{eq:QSD-int-eqn}
\lambda_A\,q_A(y)
&=
\int_0^A q_A(x)\,\mathcal{K}_{\tinyinfty}(x,y)\,dx,
\;\;\text{subject to}\;\;
\int_0^A q_A(x)\,dx=1
\end{align}
(cf.~\cite{Moustakides+etal:SS11} and~\cite{Pollak:AS85}). We note that $q_A(x)$ and $\lambda_A$ are both unique. Once $q_A(x)$ and $\lambda_A$ are found, one can compute  the ARL to false alarm $\bar{\ell}=\EV_{\tinyinfty}[\mathcal{T}_A^Q]$ and the detection delay  $\bar{\delta}=\EV_0[\mathcal{T}_A^Q]$, which is independent from the change-point:
\begin{align*}
\bar{\ell}
&=
\int_0^A\ell(x)\,q_A(x)\,dx=1/(1-\lambda_A)
\;\;\text{and}\;\;
\bar{\delta}
=
\int_0^A\delta_0(x)\,q_A(x)\,dx .
\end{align*}

The second equality in the above formula for $\bar{\ell}$ is due to the fact that, by design, the $\Pr_{\tinyinfty}$-distribution of the discrete random variable $\mathcal{T}_A^Q$ is exactly geometric with parameter $1-\lambda_A$, $0<\lambda_A<1$, i.e.,  $\Pr_{\tinyinfty}(\mathcal{T}_A^Q>\nu)=\lambda_{A}^{\nu}$, $\nu\ge0$; in general, $\lim_{A\to\infty}\lambda_A=1$. Also,~\cite{Pollak+Tartakovsky:AP11} provide sufficient conditions for $\lambda_A$ to be an increasing function of $A$; in particular, they show that if the cdf of $Z_1=\log\LR_1$ under measure $\Pr_{\tinyinfty}(\cdot)$ is concave, $\lambda_A$ is increasing in $A$.

Equations~\eqref{eq:ARL+ADD0-int-eqn}--\eqref{eq:QSD-int-eqn}  govern all of the desired performance measures (OC-s) for all considered detection procedures. The final question is that of solving these equations to evaluate the OC-s. Usually these equation cannot be solved analytically, so that  a numerical solver is  in order. We will rely on one offered by~\cite{Moustakides+etal:SS11}. Specifically, it is a piecewise-constant (zero-order polynomial) collocation method with the interval of integration $[0,A]$ partitioned into $N\ge1$ equally long subintervals. The collocation nodes are the subintervals' middle points. As the simplest case of the piecewise collocation method (see, e.g.,~\citealp[Chapter~12]{Atkinson+Han:Book09}), the question of accuracy is a well-understood one, and tight error bounds can be easily obtained from, e.g.,~\cite[Theorem~12.1.2]{Atkinson+Han:Book09}. Specifically, it can be shown that the uniform $\mathtt{L}_{\tinyinfty}$ norm of the difference between the exact solution and the approximate one is $O(1/N)$, provided $N$ is sufficiently large.

We are now in a position to apply the above performance evaluation methodology to the $\mathcal{N}(\mu,a\mu)$-to-$\mathcal{N}(\theta,a\theta)$ change-point scenario, and compute and compare against one another the performance of the four detection procedures.

\section{Performance analysis}
\label{sec:analysis} %

In this section, we  utilize the performance evaluation methodology of the preceding section for the $\mathcal{N}(\mu,a\mu)$-to-$\mathcal{N}(\theta,a\theta)$ change-point scenario~\eqref{eq:Gauss2Gauss-model} and evaluate operating characteristics of four detection procedures.  We first describe how we performed the computations and then report and discuss the obtained results.

We start with deriving pre- and post-change distributions $P_d^{\LR}(t)=\Pr_d(\LR_1\le t)$, $d=\{0,\infty\}$ of the LR.  It can be seen that
\begin{align*}
\LR_n
&=
\frac{g(X_n)}{f(X_n)}
=
\sqrt{\frac{\mu}{\theta}}\exp\left\{-\frac{\theta-\mu}{2a}\right\}
\exp\left\{\frac{\theta-\mu}{2a\theta\mu}X_n^2\right\},
\end{align*}
and
\begin{align*}
\LR_n
&\ge
\sqrt{\frac{\mu}{\theta}}\exp\left\{-\frac{\theta-\mu}{2a}\right\},
\;\;\text{if}\;\;\theta-\mu>0,
\;\;\text{and}\;\;
\LR_n
\le
\sqrt{\frac{\mu}{\theta}}\exp\left\{-\frac{\theta-\mu}{2a}\right\},
\;\;\text{if}\;\;\theta-\mu<0.
\end{align*}
Observe that to obtain the required distributions, one is to consider separately two cases $\theta-\mu>0$ and $\theta-\mu<0$.
Analytically, both cases can be treated in a similar manner, so that we present the details only in the first case.

Let $\theta-\mu>0$. For any measure $\Pr(\cdot)$ it is straightforward to see that
\begin{align*}
\begin{split}
\Pr(\LR_n\le t)
&=
\Pr\left(X_n\le\sqrt{\theta\mu\left[\left(2\log t-\log\frac{\mu}{\theta}\right)/(\theta-\mu)+1\right]}\,\right)-\\
&\qquad\qquad\qquad\qquad-
\Pr\left(X_n\le-\sqrt{\theta\mu\left[\left(2\log t-\log\frac{\mu}{\theta}\right)/(\theta-\mu)+1\right]}\,\right),\\
&\qquad\qquad\qquad\qquad\qquad\qquad\text{for}\; t\ge\sqrt{\frac{\mu}{\theta}}\exp\left\{-\frac{\theta-\mu}{2a}\right\},
\end{split}
\end{align*}
and $\Pr(\LR_n\le t)=0$ otherwise.

Hence,
\begin{align*}
\begin{split}
P_{\tinyinfty}^{\LR}(t)
&=
\Phi_{\mu,a\mu}\left(\sqrt{\theta\mu\left[\left(2\log t-\log\frac{\mu}{\theta}\right)/(\theta-\mu)+1\right]}\,\right)-\\
&\qquad\qquad\qquad\qquad
-\Phi_{\mu,a\mu}\left(-\sqrt{\theta\mu\left[\left(2\log t-\log\frac{\mu}{\theta}\right)/(\theta-\mu)+1\right]}\,\right),\\
&\qquad\qquad\qquad\qquad\qquad\qquad\text{for}\; t\ge\sqrt{\frac{\mu}{\theta}}\exp\left\{-\frac{\theta-\mu}{2a}\right\},
\end{split}
\end{align*}
and $P_{\tinyinfty}^{\LR}(t)=0$ otherwise; hereafter,
\begin{align*}
\Phi_{\mu,\sigma^2}(x)=\int_{-\infty}^x \frac{1}{\sqrt{2\pi\sigma^2}}\exp\left\{-\frac{(t-\mu)^2}{2\sigma^2}\right\}dt
\end{align*}
is the cdf of a Gaussian distribution with mean $\mu$ and variance $\sigma^2$.

Similarly,
\begin{align*}
\begin{split}
P_0^{\LR}(t)
&=
\Phi_{\theta,a\theta}\left(\sqrt{\theta\mu\left[\left(2\log t-\log\frac{\mu}{\theta}\right)/(\theta-\mu)+1\right]}\,\right)-\\
&\qquad\qquad\qquad\qquad
-\Phi_{\theta,a\theta}\left(-\sqrt{\theta\mu\left[\left(2\log t-\log\frac{\mu}{\theta}\right)/(\theta-\mu)+1\right]}\,\right),\\
&\qquad\qquad\qquad\qquad\qquad\qquad\text{for}\; t\ge\sqrt{\frac{\mu}{\theta}}\exp\left\{-\frac{\theta-\mu}{2a}\right\},
\end{split}
\end{align*}
and $P_0^{\LR}(t)=0$ otherwise.

As soon as both distributions, $P_0^{\LR}(t)$ and $P_\infty^{\LR}(t)$, are expressed in a closed form, one is ready to employ the numerical method of~\cite{Moustakides+etal:SS11}. We implemented the method in MATLAB, and solved the corresponding equations with an accuracy of a fraction of a percent.

It is also easily seen that the Kullback--Leibler information numbers $I_f=-\EV_{\tinyinfty}[\log \Lambda_1]$ and $I_g=\EV_0[\log\Lambda_1]$ are
\begin{align*}
I_f
=
\frac{(\mu-\theta)^2}{2a\theta}+\frac{1}{2}\left[\left(\frac{\mu}{\theta}-1\right)-\log\frac{\mu}{\theta}\right]
\;\;\text{and}\;\;
I_g
&=
\frac{(\theta-\mu)^2}{2a\mu}+\frac{1}{2}\left[\left(\frac{\theta}{\mu}-1\right)-\log\frac{\theta}{\mu}\right] .
\end{align*}
Also, for this model, we have
\begin{align*}
S_k
&=
\sum_{j=1}^k\log\LR_j
=
\frac{\theta-\mu}{2\mu}\sum_{j=1}^k \frac{X_j^2}{a\theta}
-\frac{k}{2}\left[\frac{\theta-\mu}{a}+\log\frac{\theta}{\mu}\right],\;\; k \ge 1 .
\end{align*}

We now consider the question of computing $\varkappa$ and $\zeta$. To this end, in general both can be computed either using~\eqref{eq:zeta+kappa}, or by Monte Carlo simulations. However, in our particular case, the latter is better since each $S_k$ is non-central chi-squared distributed, and this distribution is an infinite series, which is slowly converging for weak changes, i.e., when $\theta$ and $\mu$ are close to each other. This translates into a computational issue. Hence, $\varkappa$ and $\zeta$ as well as $C_0$ and $C_{\tinyinfty}$ are all evaluated by Monte Carlo simulations.

Specifically, to compute $\varkappa$ and $\zeta$ we generate $10^6$ trajectories of $S_k$ for $k$ ranging from 1 to $10^4$ for each measure $\Pr_d(\cdot)$, $d=\{0,\infty\}$, and estimate $\varkappa$ and $\zeta$ using~\eqref{eq:zeta+kappa}. That is, we cut the infinite series appearing in both formulae at $10^4$. The same data are used to compute constants $\beta_0$ and $\beta_{\tinyinfty}$, as well as constants $C_0$ and $C_{\tinyinfty}$. We directly estimate the corresponding integrals, which is a textbook Monte Carlo problem.

One of the most important issues is the design of the SR--$r$ procedure, namely, choosing the starting point $r$. As we discussed in~\autoref{sec:asymptotics}, the head start  $R_0^r=r^*$ should be set to the solution of the constrained minimization problem~\eqref{eq:rstar}. Clearly, with this initialization the SADD of the SR--$r$ procedure is as close to the optimum as possible. The challenge is that unlike other procedures, one cannot determine the detection threshold from the equation $\ARL(\mathcal{S}_A^r)=\gamma$ for the desired $\gamma>1$, even using the approximation $\ARL(\mathcal{S}_A^r)\approx(A/\zeta)-r$. The reason is that we now have two unknowns -- $A$ and $r$. This is where $r^*$ steps in. For a fixed $\gamma$, it is clear that as $A$ grows, in order to maintain $\ARL(\mathcal{S}_A^r)$ at the same level $\gamma$ (i.e., to keep the constraint $\ARL(\mathcal{S}_A^r)=\gamma$ satisfied), $r$ has to grow as well. The smallest value of $r$ is $0$, and when $r=0$, we can easily determine the (smallest) value of $A$ for which $\ARL(\mathcal{S}_A^r)=\gamma$, because in this case it is the SR procedure, and $\ARL(\mathcal{S}_A)\approx A/\zeta$. Once this $A$ is found, we compute $\mathcal{J}_{\mathrm{P}}(\mathcal{S}_A^r)$ and $\mathcal{J}_{\mathrm{LB}}(\mathcal{S}_A^r)$, and record the difference $\mathcal{J}_{\mathrm{P}}(\mathcal{S}_A^r)-\mathcal{J}_{\mathrm{LB}}(\mathcal{S}_A^r)$. Then we increase the detection threshold, and find $r$ for which $\ARL(\mathcal{S}_A^r)=\gamma$ is satisfied, and again compute $\mathcal{J}_{\mathrm{P}}(\mathcal{S}_A^r)$ and $\mathcal{J}_{\mathrm{LB}}(\mathcal{S}_A^r)$, and record the difference $\mathcal{J}_{\mathrm{P}}(\mathcal{S}_A^r)-\mathcal{J}_{\mathrm{LB}}(\mathcal{S}_A^r)$. This process continues until enough data is collected to figure out what $A$ and $r^*$ are. See also~\citet[Section~5]{Polunchenko+Tartakovsky:MCAP11}.

Though $r^*$ is an intuitively appealing initialization strategy, it is interesting only from a theoretical point of view. For practical purposes,~\cite{Tartakovsky+etal:TPA11} suggest to choose $r$ so as to equate $\ADD_0(\mathcal{S}_A^r)$ and $\ADD_{\tinyinfty}(\mathcal{S}_A^r)$. This is the expected behavior of the SR--$r$ procedure for large $A$'s, and we in fact also do observe it for our example. For large values of $A$, $\ADD_0(\mathcal{S}_A^r)$ can be approximated as
\begin{align*}
\ADD_0(\mathcal{S}_A^r)
&\approx
\frac{1}{I_g}(\log A+\varkappa-C_r),
\end{align*}
where
\begin{align*}
C_r
&=
\EV[\log(1+r+\tilde{V}_{\tinyinfty})]
=
\int_0^\infty\log(1+r+y)\,d\tilde{Q}(y),
\end{align*}
and $\tilde{V}_{\tinyinfty}$ and $\tilde{Q}(y)$ are as in~\autoref{sec:asymptotics}. We also recall that $\ADD_{\tinyinfty}(\mathcal{S}_A^r)=\mathcal{J}_{\rm P}(\mathcal{S}_A^r)$ admits the same asymptotics \eqref{eq:SADDSRr}, except $C_r$ is replaced with $C_{\tinyinfty}$ defined in~\eqref{eq:C-inf-def}. Hence, for large values of $A$, the equation $\ADD_0(\mathcal{S}_A^r)=\ADD_{\tinyinfty}(\mathcal{S}_A^r)$, where the unknown is $r$, reduces to the equation $C_r=C_{\tinyinfty}$. Let $r_\star$ be its solution. The important feature of $r_\star$ compared to $r^*$ is that it does not depend of the threshold (i.e., on $\gamma$). \cite{Tartakovsky+etal:TPA11} and \cite{Polunchenko+Tartakovsky:MCAP11} offer an example where $C_r$ and $C_{\tinyinfty}$ can be found in a closed form, and the equation $C_r=C_{\tinyinfty}$ can be written out and solved explicitly.

As a first illustration, suppose $\mu=1000$, $\theta=1001$, and $a=0.01$, and set $\gamma=10^4$ for each procedure in question. The corresponding detection threshold and head start $r^*$ for each procedure are reported in~\autoref{tab:A+headstart+ARL-m1000-t1001-a0.01-g10000}. For this case we have $I_f\approx 5\times 10^{-2}$, $I_g\approx5\times 10^{-2}$, $\zeta\approx0.83145$, $\varkappa\approx0.22$, $C_0\approx3.59$, $C_{\tinyinfty}\approx4.5$, $\beta_0\approx-1$, and $\beta_{\tinyinfty}\approx1.3$. We see that for the SR procedure the asymptotics $\ARL(\mathcal{S}_A)\approx A/\zeta$ is again perfect. The accuracy is great for the SRP rule as well: the asymptotics $\ARL(\mathcal{S}_A^Q)\approx A/\zeta-\mu_Q$ yields a value of $9999.5$ average observations, while the computed value is $9999.845$. For the SR--$r$ procedure, the asymptotics $\ARL(\mathcal{S}_A^r)\approx A/\zeta-r$ yields a value of $9999.57$ average observations, while the computed value is $9999.875$. For the CUSUM procedure, the approximation for the ARL \eqref{eq:ARLcusum} gives a value of $10006$, while the computed one is $10001.22$.

\begin{table}[!htb]
    \renewcommand{\tabcolsep}{5mm}
    \renewcommand{\arraystretch}{1.2}
    \centering
    \caption{Thresholds,  head start values and the ARL to false alarm for the four detection procedures for $\mu=1000$, $\theta=1001$, $a=0.01$, and $\gamma=10^4$. $\ARL(\T)$ is the actual ARL to false alarm exhibited by the procedure for the selected detection threshold and head start evaluated numerically}
    \vspace{2mm}
    \begin{tabular}{|l||*{2}{c|}|c|}
        \hline
        Procedure & $A$ & Head~Start & $\ARL(\T)=\EV_{\tinyinfty}[\T]$\\
            \hline\hline
        CUSUM & $350.75$ & any $W_0\le 1$ & $10001.223$\\
            \hline
        SR & $8314.4$ & $R_0=0$ & $10000.188$\\
            \hline
        SRP & $8392.0$ & $R_0^{Q}\propto Q_A\;\;(\mu_{Q}\approx93.699)$ & $9999.845$\\
            \hline
        SR--$r$ & $8356.0$ & $R_0^r\approx50.345$ & $9999.875$\\
            \hline
    \end{tabular}
    \label{tab:A+headstart+ARL-m1000-t1001-a0.01-g10000}
\end{table}

While not presented in this table, we now discuss the accuracy of the asymptotic approximations given in Section~\ref{sec:asymptotics}. The asymptotics for
$\mathcal{J}_{\mathrm{P}}(\mathcal{S}_A)$ given by \eqref{eq:SADDSR} yields a value of $117$, which aligns well with the observed $113$. For the CUSUM procedure we obtain that the approximation for $\mathcal{J}_{\mathrm{P}}(\mathcal{C}_A)$ (see \eqref{eq:SADDCS})  gives 101 vs. 105 observed. For $\ADD_{\tinyinfty}(\mathcal{C}_A)$ given in \eqref{eq:ADDinftyCS} we have 93 vs. 95 observed. For the approximations $\mathcal{J}_{\mathrm{P}}(\mathcal{S}_A^Q)$, $\mathcal{J}_{\mathrm{P}}(\mathcal{S}_A^r)$ and $\ADD_{\tinyinfty}(\mathcal{S}_A^r)$ (see \eqref{eq:SADDSRP} and \eqref{eq:SADDSRr}) we have 93 vs. 95 observed. We can conclude that the accuracy of asymptotic approximations is satisfactory.

\begin{figure}[!htb]
    \centering
    \includegraphics[width=0.73\textwidth]{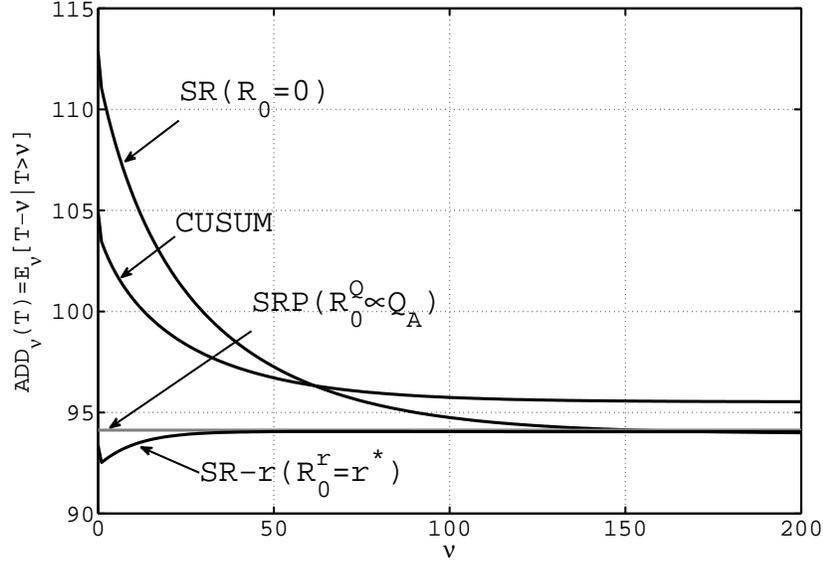}
    \vspace{-5mm}
    \caption{Conditional average detection delay  $\ADD_{\nu}(\T)=\EV_{\nu}[\T-\nu|\T>\nu]$ for the four detection procedures for $\mu=1000$, $\theta=1001$, $a=0.01$, and $\gamma=10^4$ versus the changepoint $\nu$.}
    \label{fig:ADDnu-vs-nu-m1000-t1001-a0.01-g10000}
\end{figure}

Next, the $\ADD_{\nu}(\T)=\EV_{\nu}[\T-\nu|\T>\nu]$ vs.\ the changepoint $\nu$ is shown in~\autoref{fig:ADDnu-vs-nu-m1000-t1001-a0.01-g10000}. We first discuss the procedures' performance for changes that take place either immediately (i.e., $\nu=0$), or in the near to mid-term future. As we have mentioned earlier, for the CUSUM and SR procedures the worst ADD is attained at $\nu=0$, i.e., $\mathcal{J}_{\mathrm{P}}(\mathcal{C}_A)=\EV_0[\mathcal{C}_A]$ and $\mathcal{J}_{\mathrm{P}}(\mathcal{S}_A)=\EV_0[\mathcal{S}_A]$.  From~~\autoref{fig:ADDnu-vs-nu-m1000-t1001-a0.01-g10000} we gather that at $\nu=0$ the average delay of the CUSUM procedure is $105$ observations and of the SR procedure is $113$. Overall, we see that the CUSUM procedure is superior to the SR procedure for $0\le\nu\le60$; the two procedures are at performance parity at $\nu=60$. One can therefore conclude that the CUSUM procedure is better than the SR rule for changes that occur either immediately, or in the near to mid-term future. However, for distant changes (i.e., when $\nu$ is large), it is the opposite, as expected. Below we will consider one more case where the SR procedure significantly outperforms the CUSUM procedure. For now note that the only procedure indifferent to how soon or late the change may take place is the SRP procedure $\mathcal{S}_A^Q$; recall that it is an equalizer, i.e., $\ADD_{\nu}(\mathcal{S}_A^Q)$ is the same for all $\nu\ge0$. The (flat) ADD for the SRP procedure is $94$. Since $\gamma=10^4$,  SRP and SR--$r$ are nearly indistinguishable. Yet  SR$-r$ appears to be uniformly (though slightly) better than SRP.

The obtained $\ADD_{\nu}(\T)$ for each procedure for selected values of $\nu$ are summarized in~\autoref{tab:ADDnu-vs-nu-m1000-t1001-a0.01-g10000}. The table also reports the lower bound~\eqref{eq:LowerBound-def} and the STADD $\mathcal{J}_{\mathrm{ST}}(\T)$. From the table we see that for the considered parameters all procedures deliver about the same STADD.

\begin{table}[h]
    \renewcommand{\tabcolsep}{5mm}
    \renewcommand{\arraystretch}{1.2}
    \centering
       \caption{Summary of the performance of the four procedures for $\mu=1000$, $\theta=1001$, $a=0.01$, and $\gamma=10^4$}
       \vspace{2mm}
    \begin{tabular}{|l||*{5}{c|}|c|}
        \hline
        {}& \multicolumn{5}{c||}{$\ADD_{\nu}(\T)=\EV_{\nu}[\T-\nu|\T>\nu]$} & {}\\
            \cline{2-7}
        Procedure & 0 & 50 & 100 & 150 & 200 & $\mathcal{J}_{\mathrm{ST}}(\T)$\\
            \hline\hline
        CUSUM & \begin{tabular*}{1.4cm}{@{}c}104.98\\(SADD)\end{tabular*} & 96.72 & 95.75 & 95.57 & 95.53 & 95.55\\
            \hline
        SR & \begin{tabular*}{1.4cm}{@{}c}112.87\\(SADD)\end{tabular*} & 97.26 & 94.75 & 94.15 & 94.00 & 94.00\\
            \hline
        SRP & \multicolumn{6}{c|}{94.127 (SADD)}\\
            \hline
        SR--$r$ & 93.38 & 94.04 & 94.04 & 94.04 & \begin{tabular*}{1.4cm}{@{}c}94.04\\(SADD)\end{tabular*} & 94.04 \\
            \hline\hline
        Lower~Bound & \multicolumn{6}{c|}{94.04}\\
            \hline
                \end{tabular}
    \label{tab:ADDnu-vs-nu-m1000-t1001-a0.01-g10000}
\end{table}

As another illustration, suppose $\mu=1000$, $\theta=1001$, but $a=1$. Also, assume that $\gamma=10^3$. The corresponding detection threshold and head start for each procedure are reported in~\autoref{tab:A+headstart+ARL-m1000-t1001-a1-g1000}. For this case we have $I_f\approx 5\times 10^{-4}$, $I_g\approx5\times 10^{-4}$, $\zeta\approx0.981$, $\varkappa\approx1.55$, $C_0\approx8.08$ and $C_{\tinyinfty}\approx8.82$, $\beta_0\approx-2.1$ and $\beta_{\tinyinfty}\approx2.2$. We see that for the SR procedure the asymptotics $\ARL(\mathcal{S}_A)\approx A/\zeta$ is perfect. However, the accuracy is not as great for the SRP rule: the asymptotics $\ARL(\mathcal{S}_A^Q)\approx A/\zeta-\mu_Q$ yields a value of $929.72$, while the computed value is $1000.33$. For the SR--$r$ procedure, the accuracy is slightly off as well: the asymptotics $\ARL(\mathcal{S}_A^r)\approx A/\zeta-r$ yields  $930.72$, while the computed one is $999.98$. For the CUSUM procedure, approximation \eqref{eq:ARLcusum} gives $1358.3$, while the computed one is $1000.1$.

\begin{table}[h]
    \renewcommand{\tabcolsep}{5mm}
    \renewcommand{\arraystretch}{1.2}
    \centering
    \caption{Thresholds, head start values and the ARL to false alarm  for the four detection procedures for $\mu=1000$, $\theta=1001$, $a=1$, and $\gamma=10^3$. $\ARL(\T)$ is the actual ARL to false alarm exhibited by the procedure for the selected detection threshold and head start evaluated numerically}
    \vspace{2mm}
    \begin{tabular}{|l||*{2}{c|}|c|}
        \hline
        Procedure & $A$ & Head~Start & $\ARL(\T)=\EV_{\tinyinfty}[\T]$\\
            \hline\hline
        CUSUM & $2.272$ & any $W_0\le1$ & $1000.096$\\
            \hline
        SR & $981.0$ & $R_0=0$ & $999.996$\\
            \hline
        SRP & $1844.0$ & $R_0^{Q}\propto Q_A\;\;(\mu_{Q}\approx879.248)$ & $1000.333$\\
            \hline
        SR--$r$ & $1811.0$ & $R_0^r=r^*\approx845.872$ & $999.981$\\
            \hline
    \end{tabular}
    \label{tab:A+headstart+ARL-m1000-t1001-a1-g1000}
\end{table}

The asymptotics for $\mathcal{J}_{\mathrm{P}}(\mathcal{S}_A)$ yields a value of $717$, which aligns well with the observed $722$. For the CUSUM procedure we obtain that the approximation for $\mathcal{J}_{\mathrm{P}}(\mathcal{C}_A)$ gives 541 vs. 563. For $\ADD_{\tinyinfty}(\mathcal{C}_A)$ we have 314 vs. 463. For $\mathcal{J}_{\mathrm{P}}(\mathcal{S}_A^Q)$, $\mathcal{J}_{\mathrm{P}}(\mathcal{S}_A^r)$ and $\ADD_{\tinyinfty}(\mathcal{S}_A^r)$ (which are asymptotically all the same) we have 500 vs. 502. We see that the accuracy is reasonable.

\begin{figure}[h]
    \centering
    \includegraphics[width=0.7\textwidth]{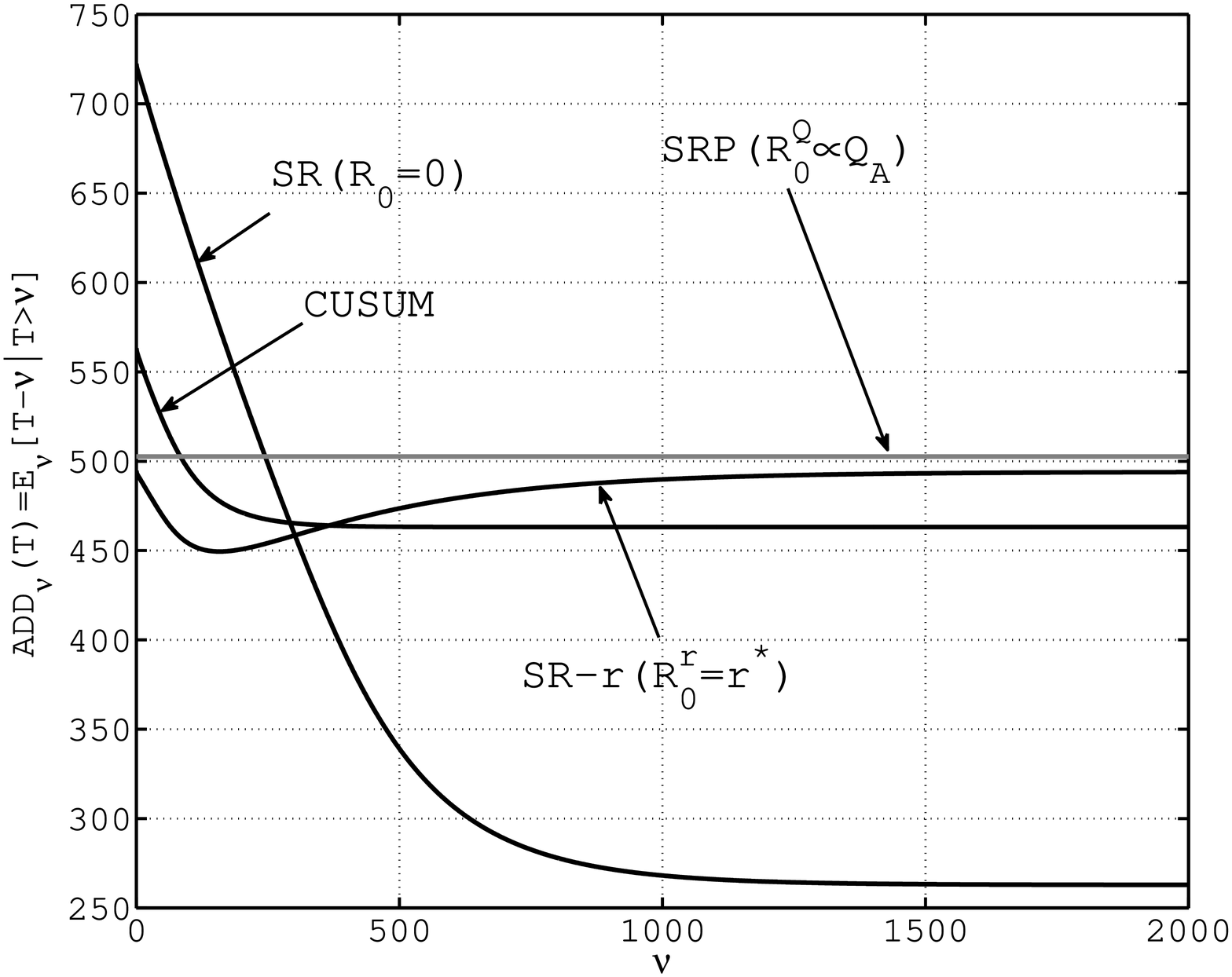}
    \vspace{-5mm}
    \caption{Average detection delay $\ADD_{\nu}(\T)=\EV_{\nu}[\T-\nu|\T>\nu]$ for the four detection  procedures for $\mu=1000$, $\theta=1001$, $a=1$, and $\gamma=10^3$.}
    \label{f:ADDnu-vs-nu-gamma=1000}
\end{figure}

 \autoref{f:ADDnu-vs-nu-gamma=1000} shows the plots of $\ADD_{\nu}(\T)=\EV_{\nu}[\T-\nu|\T>\nu]$ vs.\ the changepoint  $\nu$ for all four procedures of interest. At $\nu=0$ the CUSUM procedure requires $563$ observations (on average) to detect this change, while the SR procedure $722$ observations. Overall, we see that the CUSUM procedure is superior to the SR procedure for $0\le\nu\le270$; the two procedures are at performance parity at $\nu=270$. Again, we can conclude that the CUSUM procedure is better than the SR rule for changes that occur either immediately, or in the near to mid-term future. However, for $\nu\ge1000$, the SR procedure considerably outperforms the CUSUM procedure: the former's ADD is $263$, while the letter's is $463$. Hence, the SR procedure is better for detecting distant changes than the CUSUM procedure, as expected. The (flat) ADD for the SRP procedure is approximately $503$. This is twice as worse as the SR procedure for $\nu\ge1500$, whose detection lag for such large $\nu$'s is $263$, and stays steady for larger $\nu$'s.

Regarding the performance of the SR--$r$ procedure compared to that of the other three procedures, one immediate observation is that  SR--$r$ is uniformly better than the SRP rule; although the difference is only $8\div9$ samples, i.e., a mere fraction of a percent compared to the magnitude of the ADD itself. A similar conclusion was previously reached by~\cite{Moustakides+etal:SS11} for the Gaussian-to-Gaussian model and for the exponential-to-exponential model.
Overall, out of the four procedures in question, the best in the minimax sense is the SR--$r$ procedure. However, it is only slightly better that the SRP procedure since both are nearly asymptotical optimal in the minimax sense.

Yet another observation one may make from~\autoref{f:ADDnu-vs-nu-gamma=1000} is that the magnitude of the difference between $\ADD_0(\T)$ and $\ADD_{2000}(\T)$ for the SR procedure is extremely large, and the fact that the convergence to the steady state ADD is very slow. This can be explained by the fact that for this case the SR detection statistic is ``singular'': $R_n$ as a function of $n$ is a diagonal line with unit slope since $R_n$ deviates very little from its mean value $n$ for such small changes.

The obtained $\ADD_{\nu}(\T)$ for each procedure for selected values of $\nu$ are summarized in~\autoref{t:ADDnu-vs-nu-gamma=1000}. The table also reports the lower bound~\eqref{eq:LowerBound-def} and the STADD $\mathcal{J}_{\mathrm{ST}}(\T)$.  We remind that the SR procedure is strictly optimal in the sense of minimizing the STADD. From the table we see that its STADD is $396.44$ average observations. Since the SRP rule, $\mathcal{S}_A^Q$, is an equalizer, $\mathcal{J}_{\mathrm{P}}(\mathcal{S}_A^Q)=\mathcal{J}_{\mathrm{ST}}(\mathcal{S}_A^Q)$; for this particular case, we see that the SRP rule is as efficient as $502.64$ average observations in terms of its STADD. The other two rivals -- the CUSUM procedure and the SR--$r$ procedure -- are almost equally efficient: $471.67$ average observations for the former, and $477.43$ average observations for the latter. All in all, the SRP rule is the worst in the multi-cyclic sense.
\begin{table}[h]
    \renewcommand{\tabcolsep}{3mm}
    \renewcommand{\arraystretch}{1.2}
    \centering
    \caption{Summary of the performance of the four procedures for $\mu=1000$, $\theta=1001$, $a=1$, and $\gamma=10^3$.}
    \begin{tabular}{|l||*{7}{c|}|c|}
        \hline
        {}& \multicolumn{7}{c||}{$\ADD_{\nu}(\T)=\EV_{\nu}[\T-\nu|\T>\nu]$} & {}\\
            \cline{2-9}
        Procedure & 0 & 100 & 250 & 500 & 1000 & 1500 & 2000 & $\mathcal{J}_{\mathrm{ST}}(\T)$\\
            \hline\hline
        CUSUM & \begin{tabular*}{1.4cm}{@{}c}563.26\\(SADD)\end{tabular*} & 495.06 & 467.31 & 463.29 & 463.15 & 463.15 & 463.15 & 471.67\\
            \hline
        SR & \begin{tabular*}{1.4cm}{@{}c}722.36\\(SADD)\end{tabular*} & 626.20 & 498.64 & 339.18 & 268.14 & 263.27 & 262.91 & 396.44\\
            \hline
        SRP & \multicolumn{8}{c|}{502.636 (SADD)}\\
            \hline
        SR--$r$ & \begin{tabular*}{1.4cm}{@{}c}495.10\\(SADD)\end{tabular*} & 454.29 & 454.39 & 473.65 & 489.82 & 493.22 & 493.89 & 477.56\\
            \hline\hline
        Lower~Bound & \multicolumn{8}{c|}{485.60}\\
            \hline
    \end{tabular}
    \label{t:ADDnu-vs-nu-gamma=1000}
  \end{table}

\section{An application to cybersecurity}
\label{sec:cybersecurity}

In this section, we apply change-point detection theory in cybersecurity -- for rapid anomaly detection in computer networks' traffic.
A volume-type traffic anomaly is an event identified with a change in traffic's volume characteristics, e.g., in the traffic's intensity defined as the packet rate (i.e., number of network packets transmitted through the link per time unit). For example, a change in the packet rate may be caused by an intrusion attempt, a power outage,  or a misconfiguration in the network equipment; whether unintentionally or not, such events occur daily and are of harm to the victim. One way to mitigate the harm is by devising an automated anomaly detection system to rapidly catch suspicious, spontaneous changes in the traffic flow, so an appropriate response can be provided in a timely manner. To this end, such an anomaly detection system can be designed, e.g., by using change-point detection.

This is not a new idea; see, e.g.,~\cite{Blazek+etal:WSMC01,Tartakovsky+etal:IEEE-ToSP06,Tartakovsky+etal:SM06}. It is typically a modification of either the CUSUM procedure or Wald's Sequential Probability Ratio Test (SPRT). The SR procedure has never been paid much attention in this application. However, recall that it is exactly optimal in the multi-cyclic sense for detecting changes that occur in a distant future. We believe that this formulation is appropriate for the problem of anomaly detection in computer networks, and hence, the SR procedure may be a better option than CUSUM. We will show that this is indeed the case.

First, we note that in order to employ a detection procedure it is desirable to know the pre- and post-change models. As we focus on volume-type anomalies, the observations in our case are instantaneous values of the packet rate (traffic intensity). It is currently believed that the behavior of the packet rate can be well modeled by a Poisson process. The main reason is the tremendous growth of the Internet backbone in recent years. See, e.g.,~\cite{Cao+etal:NEC02},~\cite{Karagiannis+etal:IEEE-INFOCOM04}, and~\cite{Vishwanathy+etal:IEEE-ANTS09}. Using real data, we will show that the $\mathcal{N}(\mu,a\mu)$-to-$\mathcal{N}(\theta,a\theta)$ model (which is the standard continuous approximation to the discrete Poisson process when $a=1$) is appropriate to characterize the traffic intensity. Moreover, this model is more general than the Poisson model since by controlling parameter $a$ one can change the effect of a shift in the mean on the variance. 

We now present the results of testing of SR and CUSUM for a real Distributed Denial-of-Service (DDoS) attack, namely, an Internet Control Message Protocol (ICMP) reflector attack. The essence of this kind of attacks is to congest the victim's link with echo reply (ping) requests sent by a large number of separate compromised machines (reflectors) so as to have the victim's machine exhaust all of its resources handling the ping requests and ultimately crush. This kind of an attack clearly creates a volume-type of anomaly in the victim's traffic flow.

The ICMP data trace is courtesy of the Los Angeles Network Data Exchange and Repository (LANDER) project (see~\href{http://www.isi.edu/ant/lander/}{http://www.isi.edu/ant/lander}). This is a research-oriented traffic capture, storage and analysis infrastructure that operates under Los Nettos, a regional Internet Service Provider (ISP) in the Los Angeles area. The aggregate traffic volume through Los Nettos measures over 1~Gigabit/s each way, and the ISP's backbone is 10~Gigabit. Leveraged by a Global Positioning System (GPS) clock, LANDER's capture equipment is able to collect data at line speed and with a down-to-the-nanosecond time precision.

\begin{figure}[!htb]
    \centering
    \subfigure[Packet rate.]{\label{fig:ICMP:packet-rate}
        \includegraphics[width=0.4\textwidth]{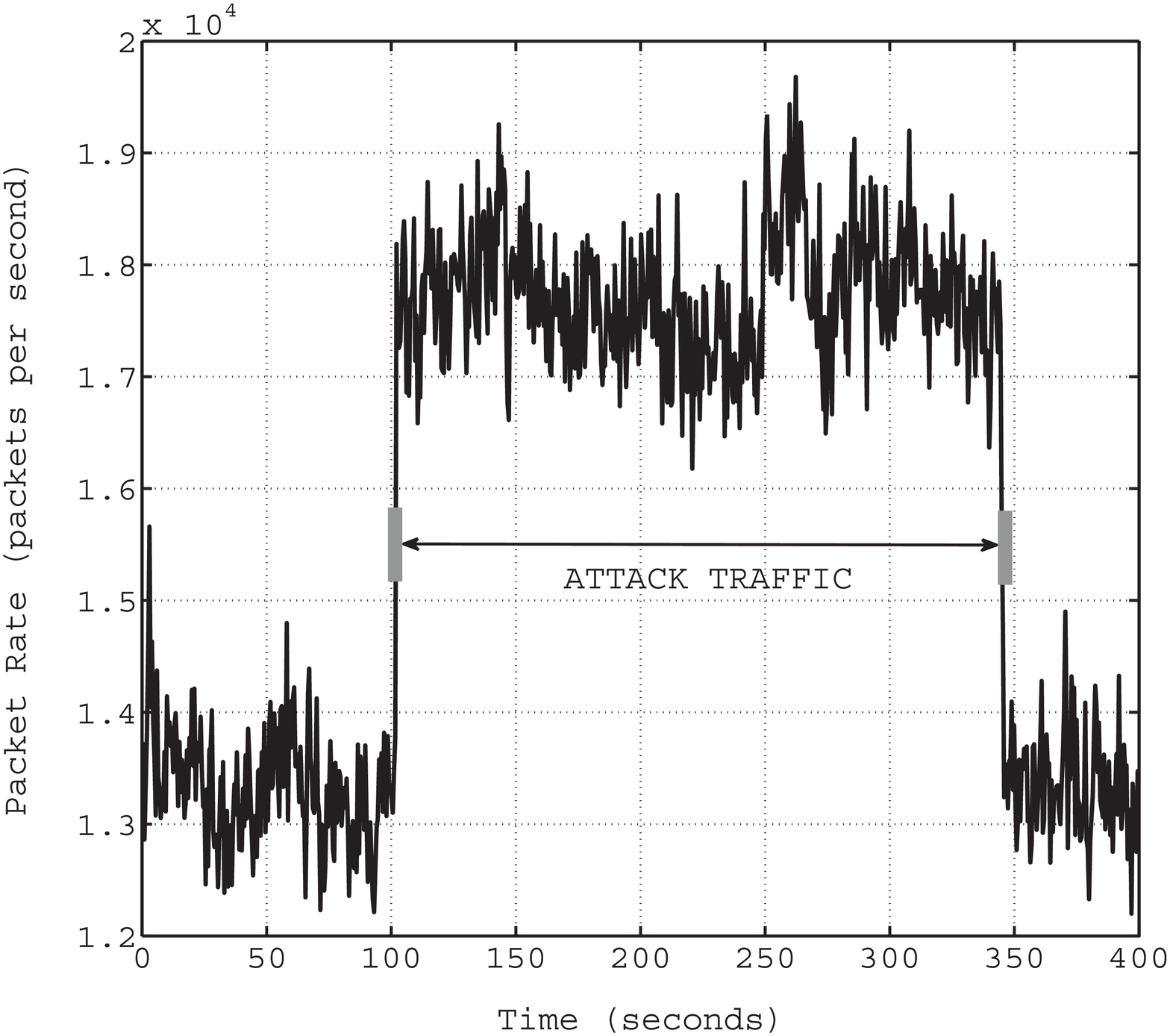}
    }
    \subfigure[Bit rate.]{\label{fig:ICMP:bit-rate}
        \includegraphics[width=0.4\textwidth]{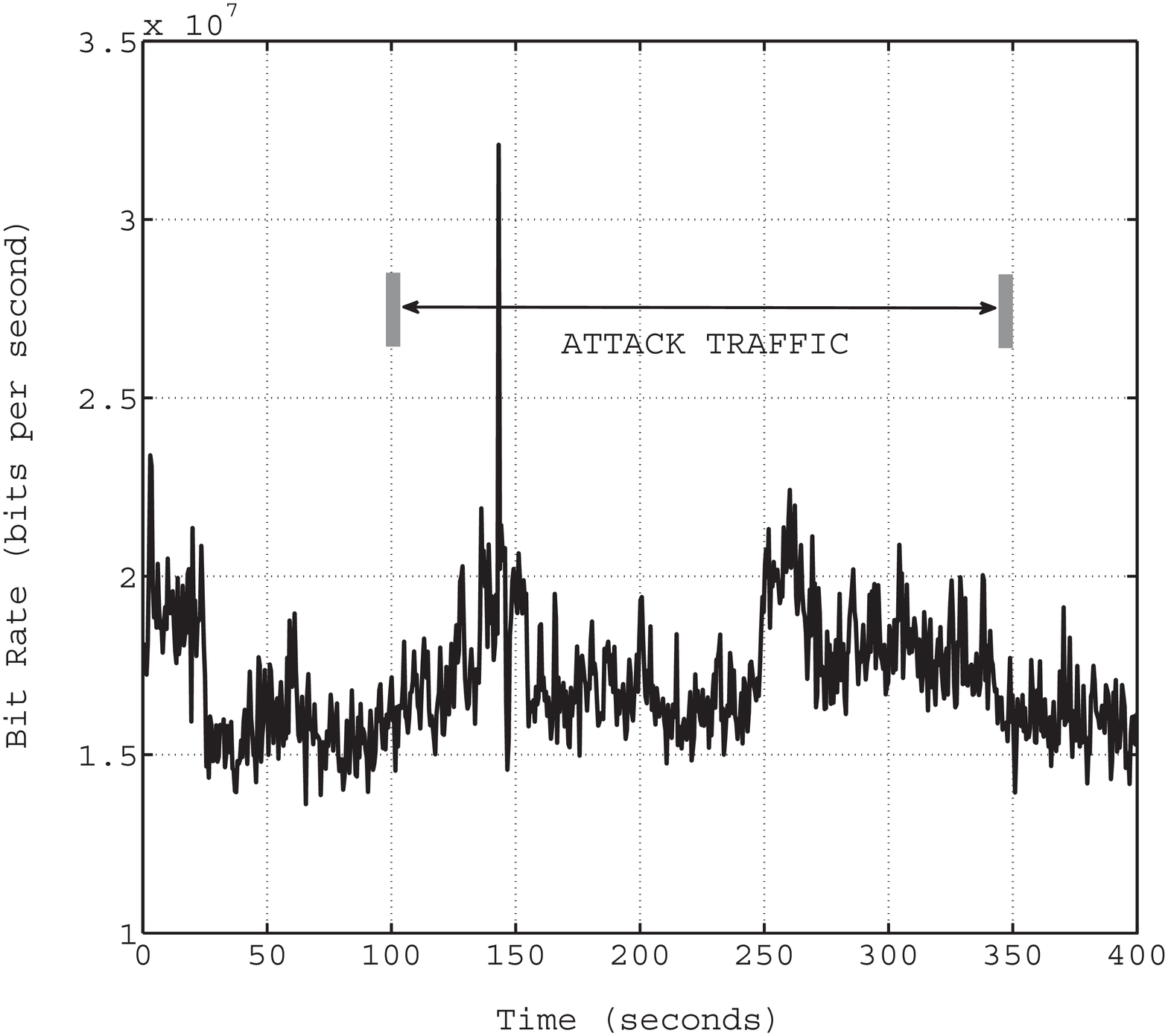}
    }
    \caption{ICMP reflector attack: traffic volume characteristics.}
    \label{fig:ICMP}
\end{figure}

The attack starts at roughly $102$ seconds into the trace and lasts for about $240$ seconds. According to the data provider, the number of reflectors is $143$, which gives a general idea as to the attack's intensity. The exact intensity can be gathered from~\autoref{fig:ICMP}, which depicts the corresponding traffic volume characteristics. Specifically, Figure~\ref{fig:ICMP:packet-rate} shows the packet rate as a function of time. The sampling rate is $0.5$ seconds. The total number of samples in the trace is $879$. It can be seen that the packet rate explodes at the attack's starting point ($102$ seconds into the trace), and goes back to the original (pre-attack) level at the attack's ending point (roughly $348$ seconds into the trace). Figure~\ref{fig:ICMP:bit-rate} shows the bit rate as a function of time. The bit rate is defined as the number of bits of information transmitted through the link per time unit. This is yet another volume characteristic. Apart from the spike at approximately $148$ seconds into the trace, its behavior is relatively calm throughout the trace. This is consistent with the fact that this is an ICMP reflector attack, i.e., a flooding-type of attack consisting of ping requests, which are small in size. That is, during the attack the traffic flow is very packet-dense, but the aggregate amount of information carried by the packets stays about the same as it was prior to the attack.

We estimated the packet rate's average and variance for legitimate traffic and for attack traffic too; in both cases, to estimate the average we used the usual sample mean, and to estimate the variance we used the usual sample variance. For legitimate traffic we obtained that the average packet rate is about $13329.764$ packets per second, and the variance is in the neighborhood of $266972.736$ packets per second; the aggregate number of samples is $401$. For attack traffic the numbers are $17723.833$ and $407968.14$, respectively, while the number of samples is $478$. We can now see the effect of the attack: it leads to an increase in the packet rate's average and variance by approximately 50\% each. Also, note that both averages are rather high. This justifies our large values for $\mu$ and $\theta$ chosen in the preceding section. This also justifies the point made in the introduction that if the packet rate truly were Poisson, with such high average it would be practically indistinguishable from Gaussian with $\sigma^2=\mu$.

\begin{figure}[!htb]
    \centering
    \subfigure[Legitimate (pre-attack) traffic.]{
        \includegraphics[width=0.4\textwidth]{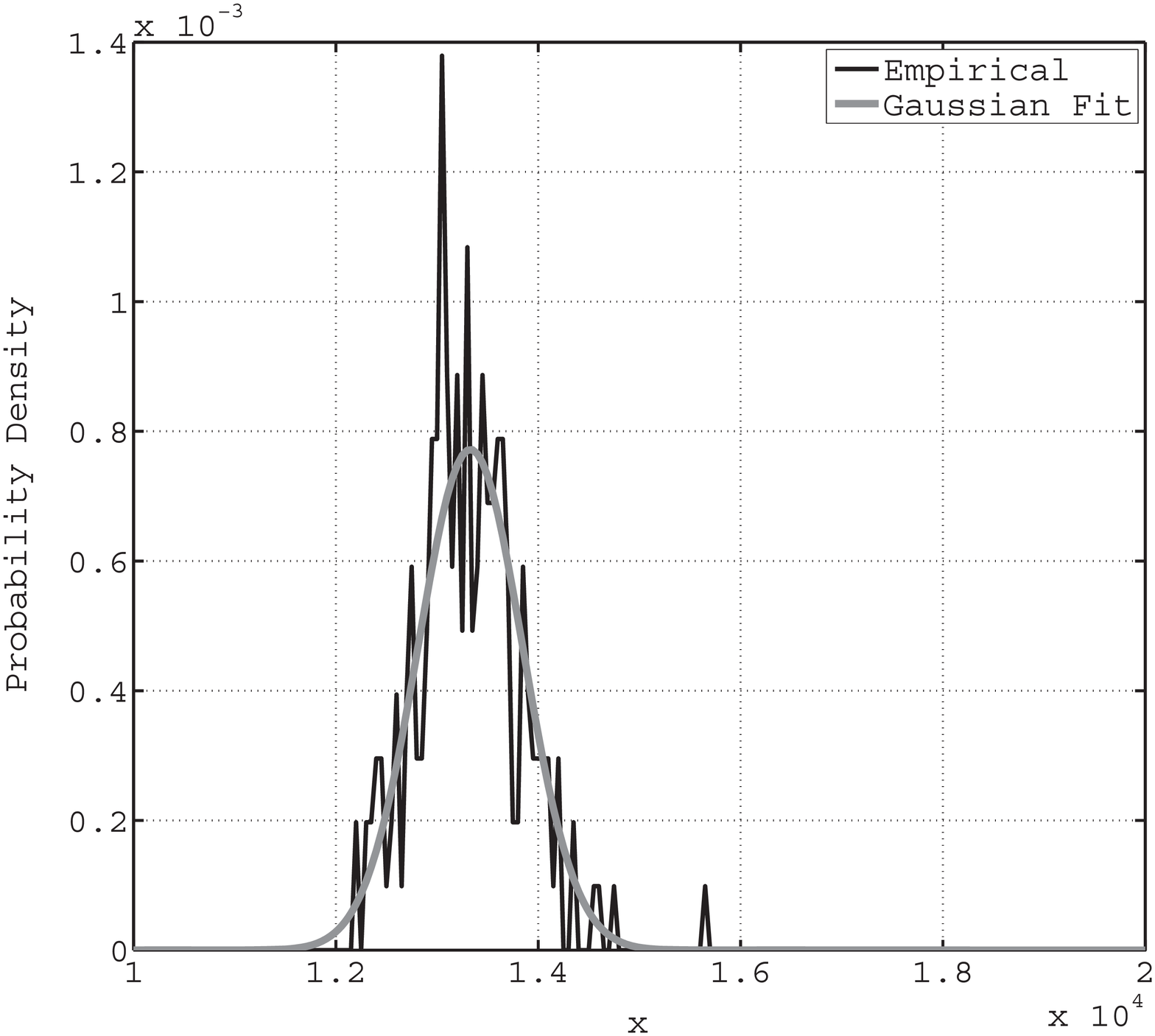}
        \label{fig:ICMP-pre-change-pdf}
    }
    \subfigure[Attack traffic.]{
        \includegraphics[width=0.4\textwidth]{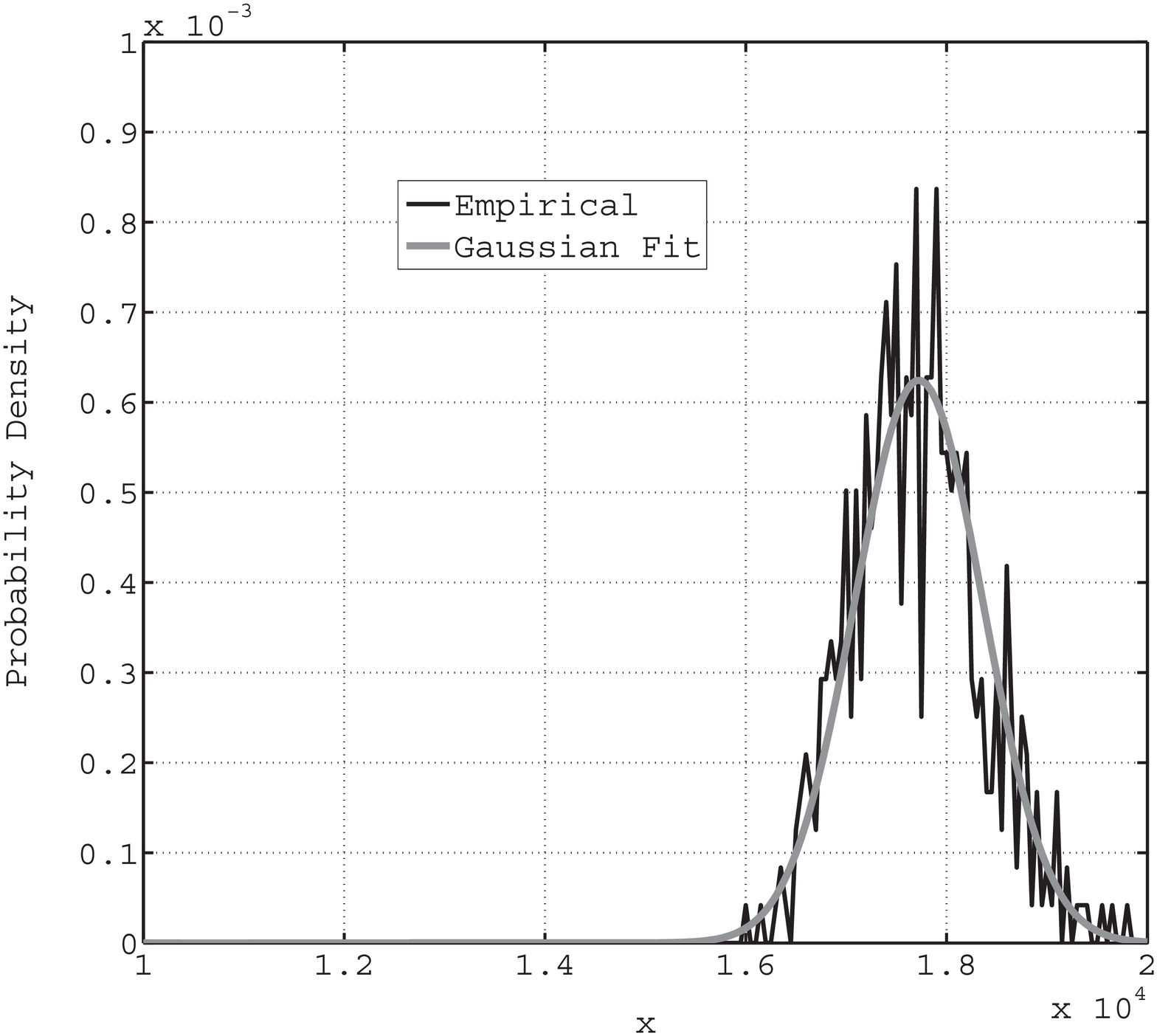}
        \label{fig:ICMP-post-change-pdf}
    }
    \caption{ICMP reflector attack: packet rate pdf with a Gaussian fit.}
    \label{fig:ICMP-pdf}
\end{figure}

The fact that neither for legitimate traffic nor for attack traffic the variance is not the same as the mean suggests to revise the conventional Poisson model. A Gaussian distribution may be a better alternative. To see how well the data agree with this hypothesis, we now access the goodness of fit of a Gaussian distribution and the data. \autoref{fig:ICMP-pdf} shows the empirical densities of the packet rate for legitimate and attack traffic. A high level of resemblance of both curves with the Gaussian distribution is apparent.

\begin{figure}[!htb]
    \centering
    \subfigure[Legitimate (pre-attack) traffic.]{
        \includegraphics[width=0.4\textwidth]{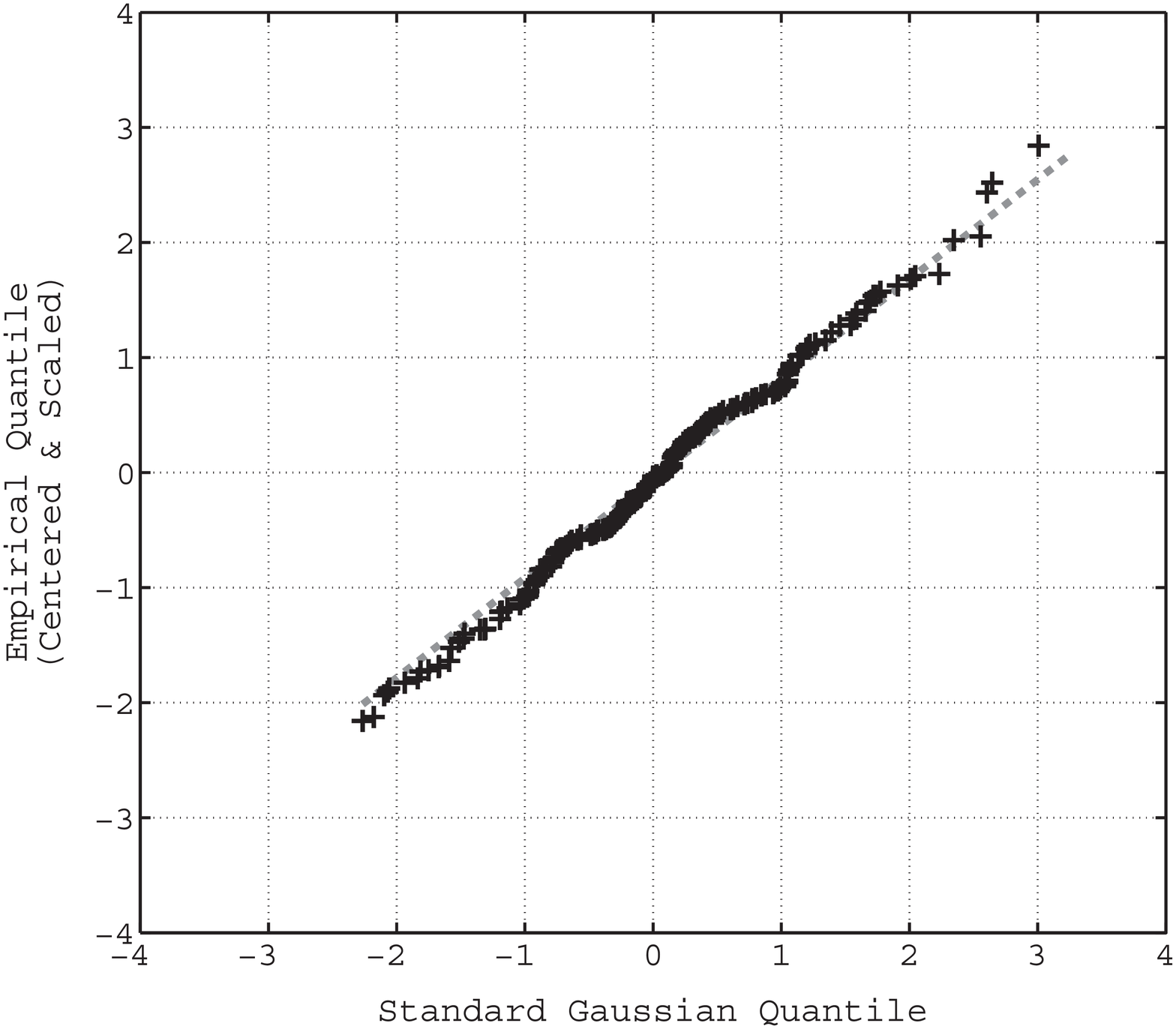}
        \label{fig:ICMP-pre-change-qq}
    }
    \subfigure[Attack traffic.]{
        \includegraphics[width=0.4\textwidth]{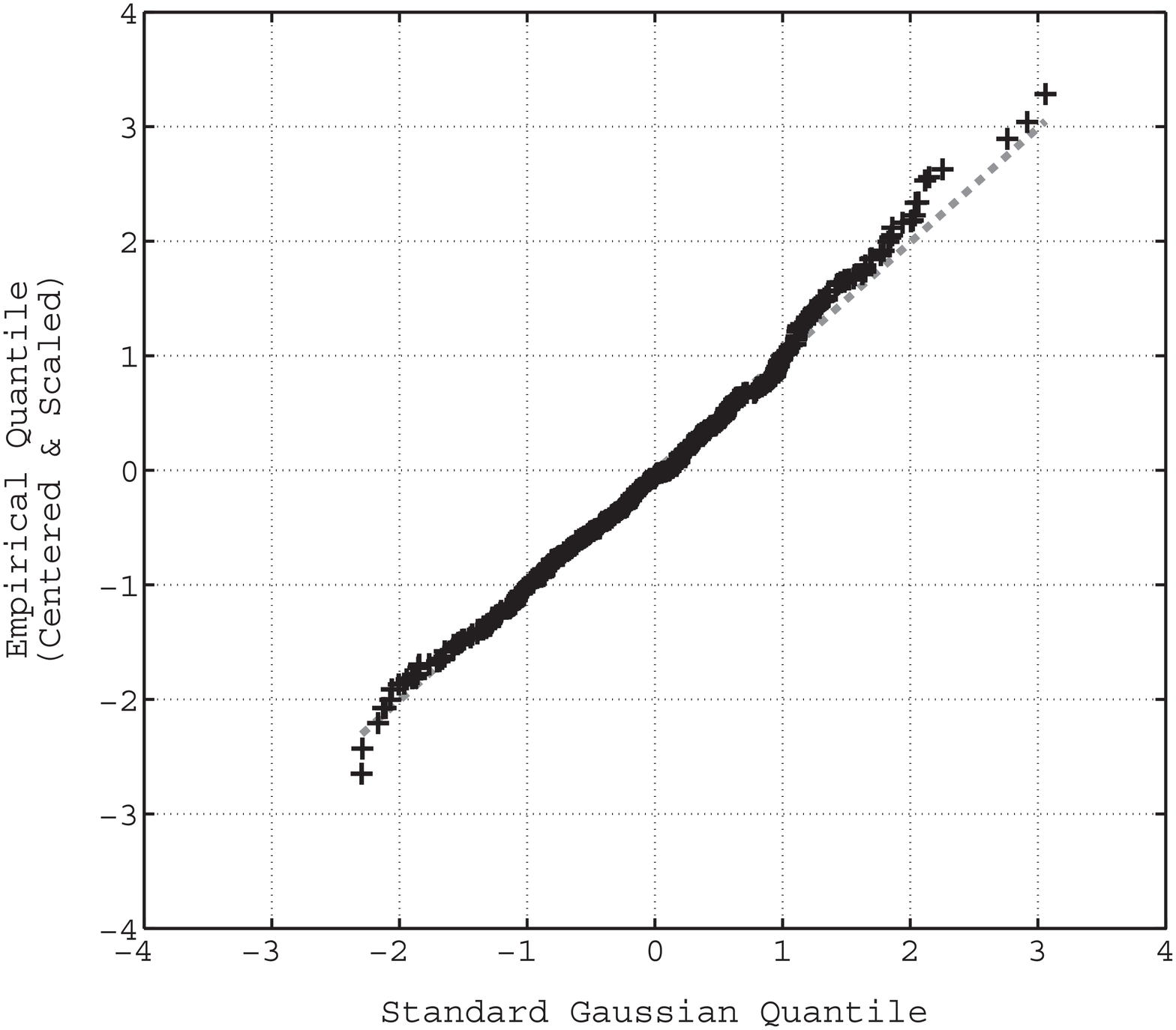}
        \label{fig:ICMP-post-change-qq}
    }
    \caption{ICMP reflector attack: Q-Q plots for packet rate distribution vs. Gaussian distribution.}
    \label{fig:ICMP-QQ}
\end{figure}

The same conclusion can be made from an inspection of the corresponding Q-Q plots (quantile-quantile) shown in~\autoref{fig:ICMP-QQ}. Specifically, the Q-Q plot for legitimate traffic is shown in Figure~\ref{fig:ICMP-pre-change-qq}, and one for attack traffic in Figure~\ref{fig:ICMP-post-change-qq}. Both plots are for centered and scaled data, so the fitted Gaussian distribution is the standard normal distribution. The fact that either plot is a straight line also confirms the ``Gaussianness'' of the data distribution.

In summary, we conclude that distributions of the packet rate are Gaussian for both legitimate and attack traffic. For legitimate traffic the variance-to-mean ratio is about $20.028$ and for attack traffic it is $23.018$. This is evidence that the $\mathcal{N}(\mu,a\mu)$-to-$\mathcal{N}(\theta,a\theta)$ model with $a=20\div23$ is appropriate to use for describing the behavior of the traffic flow in the trace.

\begin{figure}[!htb]
    \centering
    \includegraphics[width=0.7\textwidth]{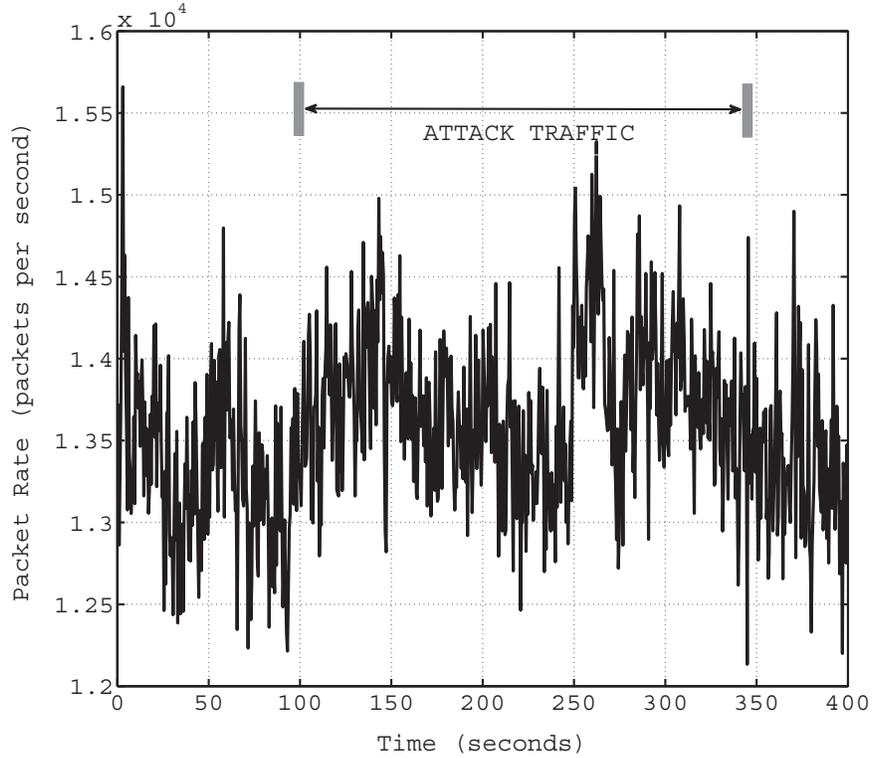}
    \caption{ICMP reflector attack: packet rate with attack's intensity diminished.}
    \label{fig:ICMP-dim-packet-rate}
\end{figure}

We now proceed to the performance analysis.  Since the ICMP attack is extremely contrast, we expect that any reasonable detection procedure will detect it in no time. As a result, it will be difficult to judge whether SR is better than CUSUM or not. This can be overcome by ``manually'' lowering the intensity of the attack, e.g., by applying the following transformation to values of the packet rate for attack's traffic
\begin{align*}
\tilde{X}_i
&=
\sqrt{13600\times 20.028}\times\cfrac{X_i-17723.833}{\sqrt{407968.14}}+13600,
\end{align*}
where $X_i$ is the original value of the attack's packet rate at the $i$-th sample. This makes the data in the trace behave according to the $\mathcal{N}(\mu,a\mu)$-to-$\mathcal{N}(\theta,a\theta)$ model with $\mu=13329.764$, $\theta=13600$ and $a=20.028$. The result is shown in~\autoref{fig:ICMP-dim-packet-rate}. We see that the attack is not as contrast any more.

\begin{figure}[!htb]
    \centering
    \subfigure[By the Shiryaev--Roberts (SR) procedure.]{
        \includegraphics[width=0.4\textwidth]{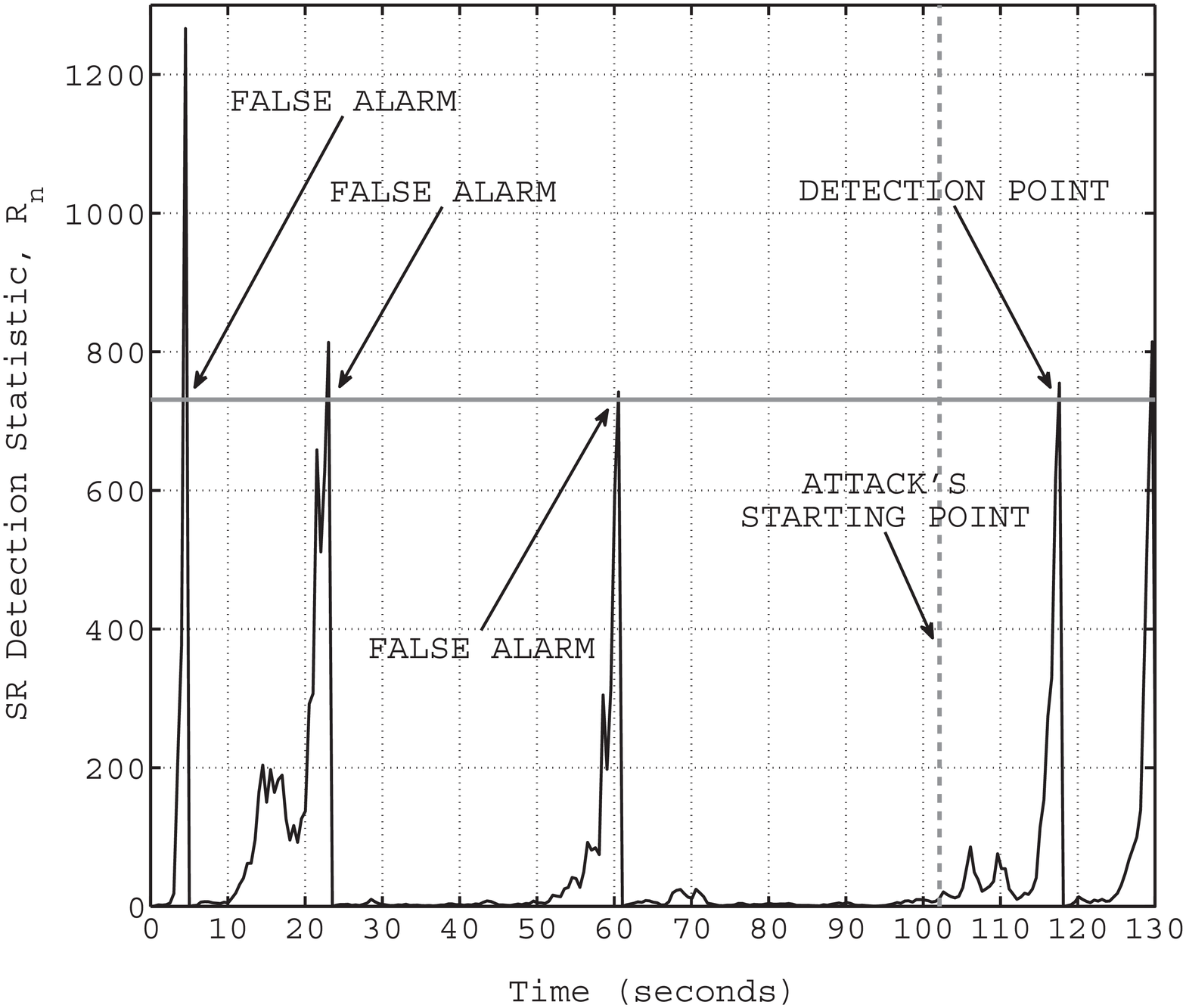}
        \label{fig:ICMP-dim-detection:SR}
    }
    \subfigure[By the CUSUM procedure.]{
        \includegraphics[width=0.4\textwidth]{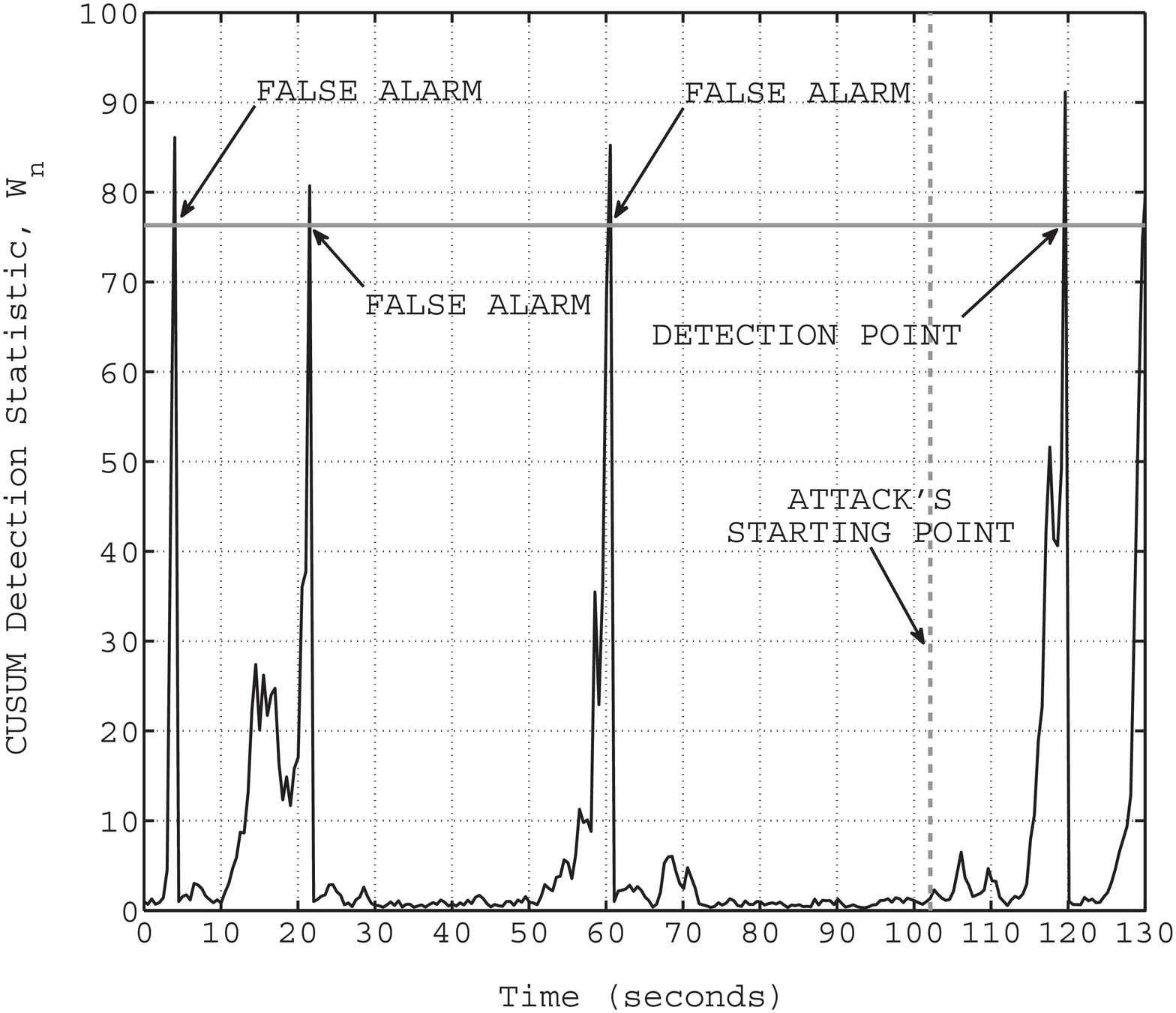}
        \label{fig:ICMP-dim-detection:CUSUM}
    }
    \caption{ICMP (diminished) reflector attack: detection.}
    \label{fig:ICMP-dim-detection}
\end{figure}

The next step is to make sure that the ARL to false alarm for each procedure is at the {\em same} desired level $\gamma>1$. To accomplish this, we recall that asymptotically, for $\gamma$ sufficiently large, $\ARL(\mathcal{C}_A)\approx A/(I_g\zeta^2)-\log A/I_f-1/(I_g\zeta)$ and $\ARL(\mathcal{S}_A)\approx A/\zeta$; in our case, $I_g\approx 0.1369$, $I_f\approx 0.1342$ and $\zeta\approx0.7313$. Thus, for, e.g., $\gamma=1000$ (which is moderate from a practical point of view), the threshold for the CUSUM procedure should be set to $76.32$, and that for the SR procedure to $731.3$. We confirmed both thresholds numerically using the methodology of~\autoref{s:perf-eval}: the actual ARL to false alarm for the CUSUM procedure is $998.4$, and that for the SR procedure is $1000.1$. We can now employ both procedures to detect the diminished version of the attack. The result is shown in~\autoref{fig:ICMP-dim-detection}. Specifically, Figure~\ref{fig:ICMP-dim-detection:SR} shows the behavior of the SR detection statistic for the first $120$ seconds of the trace's data (recall that the attack starts $102$ seconds into the trace), and Figure~\ref{fig:ICMP-dim-detection:CUSUM} shows the same for the CUSUM detection statistic.

We see that both procedures successfully detect the attack, though at the expense of raising three false alarms along the way. The detection delay for the SR procedure is roughly $15$ seconds (or $30$ samples), and that for the CUSUM procedure is about $18$ seconds (or $36$ samples). Thus, the SR procedure is better. Yet one may argue that the SR procedure is only slightly better. The difference is small simply because the trace is too short: if the attack took place at a point farther into the trace, and we had enough ``room'' to repeat the SR procedure at least a few times, the detection delay would be much smaller than that of the CUSUM procedure. Since in real life legitimate traffic dominates, the multi-cyclic idea is a natural fit.

\section{Conclusion}
\label{sec:conclusion}

We considered the basic iid version of the change-point detection problem for the Gaussian model with mean $\mu>0$ and variance $\sigma^2$ connected  via $\sigma^2=a\mu$ with $a>0$, a known constant. Of the two degrees of freedom -- $\mu$ and $a$ -- the change is only in the mean (from one known value to another), though the variance is affected as well. For this scenario, we accomplished two objectives.

First, we carried out a peer-to-peer multiple-measure comparative performance analysis of four detection procedures: CUSUM, SR, SRP and SR--$r$. We benchmarked the performance of the procedures  via Pollak's maximal average delay to detection and Shiryaev's multi-cyclic stationary average delay to detection,  subject to a tolerable lower bound, $\gamma$, on the ARL to false alarm. We tested the asymptotic (as $\gamma\to\infty$) approximations for operating characteristics using numerical techniques (solving Fredholm-type integral equations numerically). The overall conclusion is that in practically interesting cases the approximations' accuracy is ``reasonable'' or better.

Second, we considered an application of change-point detection to cybersecurity, specifically to the problem of rapid anomaly detection in computer networks' traffic. Using real data we first showed that network traffic's statistical behavior can be well-described by the Gaussian model with $\sigma^2=a\mu$; due to the fact that the effect of a shift in the mean on the variance can be controlled through parameter $a$, the Gaussian model is more flexible than the traditional Poisson one. We  then successively employed the SR and CUSUM procedures to detect an ICMP reflector attack. The SR procedure was quicker in detecting the attack, which can be explained by its  exact multi-cyclic optimality when detecting changes occurring in a far time horizon. Based on this we recommend the SR procedure for the purposes of anomaly detection in computer networks.

\section*{Acknowledgements}

The work of Aleksey Polunchenko and Alexander Tartakovsky was supported by the U.S.\ Army Research Office under MURI grant  W911NF-06-1-0044, by the U.S.\ Air Force Office of Scientific Research under MURI grant FA9550-10-1-0569, by the U.S.\ Defense Threat Reduction Agency under grant HDTRA1-10-1-0086, and by the U.S.\ National Science Foundation under grants CCF-0830419 and EFRI-1025043 at the University of Southern California, Department of Mathematics. We also thank Christos Papadopoulos of Colorado State University and John Heidemann of the Information Sciences Institute for providing real computer network data traces.


\end{document}